\documentclass[10pt]{article}
\usepackage{graphicx,float,latexsym,amsmath,amsfonts,amssymb,subeqnarray,colordvi,epstopdf,bbm}
\usepackage[latin1]{inputenc}
\usepackage{a4wide}
\usepackage{tikz}
\usepackage[hidelinks]{hyperref}

\def\R{\mathbb{R}}

\newtheorem{theorem}{Theorem}

\title{An adjoint method for the exact calibration\\ of Stochastic Local Volatility models}
\author{Maarten Wyns\footnote{Department of Mathematics and Computer Science,
University of Antwerp, Middelheimlaan 1, B-2020 Antwerp, Belgium.
\mbox{Email}: \texttt{\{maarten.wyns,karel.inthout\}@uantwerpen.be}.}
~and Karel~J.~in 't Hout\footnotemark[\value{footnote}]
}
\date{\today}

\begin{document}

\maketitle

\begin{abstract}
\noindent
This paper deals with the exact calibration of semidiscretized stochastic local volatility (SLV) models to their underlying semidiscretized local volatility (LV) models. 
Under an SLV model, it is common to approximate the fair value of European-style options by semidiscretizing the backward Kolmogorov equation using finite differences. 
In the present paper we introduce an adjoint semidiscretization of the corresponding forward Kolmogorov equation. 
This adjoint semidiscretization is used to obtain an expression for the leverage function in the pertinent SLV model such that the approximated fair values defined by the LV and SLV models are identical for non-path-dependent European-style options. 
In order to employ this expression, a large non-linear system of ODEs needs to be solved. The actual numerical calibration is performed by combining ADI time stepping with an inner iteration to handle the non-linearity. 
Ample numerical experiments are presented that illustrate the effectiveness of the calibration procedure.
\end{abstract}
\vspace{0.2cm}\noindent
\vspace{3mm}
\normalsize

\setcounter{equation}{0}
\section{Introduction}\label{intro}

In contemporary financial mathematics, \textit{stochastic local volatility} (SLV) \textit{models} are \mbox{state-of-the-art} for describing asset price processes, notably foreign exchange (FX) rates, see e.g.\ \cite{L02,TF10}. They constitute a natural combination of \textit{local volatility} (LV) and \textit{stochastic volatility} (SV) models. 
Denote by $S_{\tau}>0$ the FX rate at time $\tau \ge 0$ and consider the standard transformed variable $X_{\tau} = \log(S_{\tau}/S_{0})$. 
We deal in this paper with general SLV models of the type
\begin{equation}
\left\{ \begin{array}{l}
dX_{\tau} = (r_{d} - r_{f} - \tfrac{1}{2}\sigma^{2}_{SLV} (X_{\tau},\tau) \psi^{2}(V_{\tau})) d\tau + \sigma_{SLV} (X_{\tau},\tau) \psi(V_{\tau}) dW^{(1)}_{\tau}, \\\\
dV_{\tau} = \kappa (\eta - V_{\tau}) d\tau + \xi V_{\tau}^{\alpha} dW^{(2)}_{\tau},
\end{array} \right.
\label{eq:SLVmodelX}
\end{equation}
with $\psi(v)$ a non-negative function, $\alpha$ a non-negative parameter, $\kappa, \eta, \xi$ strictly positive parameters, $dW^{(1)}_{\tau} \cdot dW^{(2)}_{\tau} = \rho d\tau$, $-1 \leq \rho \leq 1$, and given spot values $S_{0}$, $V_{0}$. The function $\sigma_{SLV} (x,\tau)$ is often called \textit{the leverage function} and $r_{d}$, respectively $r_{f}$, denotes the risk-free interest rate in the domestic currency, respectively foreign currency. 
The SLV model \eqref{eq:SLVmodelX} can be viewed as obtained from a mixture of the LV model
\begin{equation}
dX_{LV,\tau} = (r_{d} - r_{f}-\tfrac{1}{2}\sigma^{2}_{LV} (X_{LV,\tau},\tau)) d\tau + \sigma_{LV} (X_{LV,\tau},\tau)dW_{\tau},
\label{eq:LVmodelX}
\end{equation}
with LV function $\sigma_{LV}(x,\tau)$, and the SV model 
\begin{equation}
\left\{ \begin{array}{l}
dX_{SV,\tau} = (r_{d} - r_{f} - \psi^{2}(V_{SV,\tau})) d\tau + \psi(V_{SV,\tau}) dW^{(1)}_{\tau}, \\\\
dV_{SV,\tau} = \kappa (\eta - V_{SV,\tau}) d\tau + \xi_{SV} V_{SV,\tau}^{\alpha} dW^{(2)}_{\tau}.
\end{array} \right.
\label{eq:SVmodelX}
\end{equation}
Clearly, if $\sigma_{SLV} (x,\tau)$ is identically equal to one, then the SLV model reduces to a SV model. Next, if the stochastic volatility parameter $\xi$ is equal to zero, then the SLV model reduces to~a~LV~model. 

The choice $\psi(v)=\sqrt{v}$, $\alpha = 1/2$ corresponds to the well-known Heston-based S(L)V model, the choice $\psi(v)=v$, $\alpha = 1$ to the S(L)V model considered in \cite{TF10} and the choice $\psi(v)=\exp(v)$, $\alpha = 0$ corresponds to the S(L)V model based on the exponential Ornstein--Uhlenbeck model described in \cite{S87}.

If $\alpha$ is strictly positive, we assume that $\psi(0)=0$ and the processes $V_{\tau}, V_{SV,\tau}$ are non-negative.
For $0 < \alpha < 1/2$ it holds that $V_{\tau}=0$ is attainable, for $\alpha > 1/2$ it holds that $V_{\tau}=0$ is unattainable, and for $\alpha = 1/2$ one has that $V_{\tau}=0$ is attainable if $2\kappa \eta < \xi^{2} $, see e.g.\ \cite{AP07}. The analogous result is true for the pure SV model \eqref{eq:SVmodelX}.

In financial practice, $\sigma_{LV}(x,\tau)$ is determined such that the LV model \eqref{eq:LVmodelX} yields the exact market prices for vanilla options, see e.g.\ \cite{AH11, D94}, and the parameters $\kappa, \eta, \xi_{SV}$ are chosen such that the SV model \eqref{eq:SVmodelX} reflects the market dynamics of the underlying asset, see e.g.\ \cite{TF10}. 
For the calibration of the SLV model \eqref{eq:SLVmodelX} it is customary to start from the parameters of the SV model under consideration, and to define $\xi = \mu \xi_{SV}$ with mixing parameter \mbox{$0 \le \mu \le 1$}. The mixing parameter tunes between the local volatility and stochastic volatility features. 
Next, the leverage function $\sigma_{SLV}$ is calibrated such that the SLV model yields the exact market prices for European call and put options. In the literature, no closed-form analytical relationship appears to be available between the leverage function and the fair value of vanilla options within the SLV model.
Accordingly, in financial practice the leverage function is calibrated by making use of a relationship between the SLV model and the LV model. 
It is well-known, see e.g.\ \cite{G86,T11}, that these models yield the same marginal distribution for the exchange rate $S_{\tau}$, and hence always define the same fair value for vanilla options, if the leverage function $\sigma_{SLV}(x,\tau)$ satisfies
\begin{equation}
\sigma^{2}_{LV} (x,\tau) = \mathbb{E}[\sigma_{SLV}^{2} (X_{\tau},\tau)\psi^{2}(V_{\tau}) \vert X_{\tau} = x] = \sigma_{SLV}^{2} (x,\tau)\mathbb{E}[\psi^{2}(V_{\tau}) \vert X_{\tau} = x],
\label{eq:SigmatoMatch}
\end{equation}
for all $x \in \R, \tau \ge 0$.
The latter conditional expectation can be written as
\begin{equation} 
\mathbb{E}[\psi^{2}(V_{\tau}) \vert X_{\tau} = x] = \frac{\int_{-\infty}^{\infty} \psi^{2}(v) p(x,v,\tau;X_{0},V_{0})dv}{\int_{-\infty}^{\infty}p(x,v,\tau;X_{0},V_{0})dv}, 
\label{eq:CondExpec}
\end{equation} 
where $p(x,v,\tau;X_{0},V_{0})$ denotes the joint density of $(X_{\tau},V_{\tau})$ given by the SLV model.
Since the LV function is determined such that the LV model yields exactly the observed market prices for vanilla options, the SLV model will also exactly define the same fair value whenever one is able to determine the conditional expectation above and defines the leverage function by \eqref{eq:SigmatoMatch}.
This conditional expectation itself depends on $\sigma_{SLV}(x,\tau)$, however, and determining it is a highly non-trivial task.
Recently, a variety of numerical techniques, see e.g.\ \cite{C11,EKO12,H09,RMQ07,VGO14}, has been proposed in order to approximate this conditional expectation and to approximate the appropriate leverage function. 

The numerical techniques presented in the references above do not take into account explicitly that, even if the LV function is known analytically, it is often not possible to determine exactly the corresponding fair value of vanilla options. Even within the LV model one relies on numerical methods in order to approximate the fair option values. A common approach consists of numerically solving the corresponding backward partial differential equation (PDE) by for example finite difference or finite volume methods, see e.g.\ \cite{TR00}.
When calibrating the SLV model to the LV model,
the best result one can thus aim for is that the numerical approximation of the fair value of vanilla options is the same for both models whenever similar numerical valuation methods are used.

In this paper we assume that the fair option value (within the LV model, resp. SLV model) is approximated through numerically solving the backward PDE (corresponding to the LV model, resp. SLV model) by standard finite difference or finite volume methods. Given such a spatial discretization for the backward PDE, an \textit{adjoint spatial discretization}\, will be introduced for the corresponding forward PDE. 
This adjoint spatial discretization has the important property that it always defines \textit{exactly the same approximation}\, for the fair value of non-path-dependent European options as the approximation given by the discretization of the backward equation.  
Moreover, if similar spatial discretizations are used for the backward PDE associated with the LV model and the backward PDE associated with the SLV model, then their adjoint spatial discretizations can be employed to~\textit{create an exact match between the approximations for the fair value of vanilla options within the LV model and the SLV model.}

The main contributions of this article can be visualized in the following scheme:
\begin{center}
\begin{tikzpicture} 
  \path (0,0) node[align=center,draw,rounded corners](BLV) {Discretization\\ Backward PDE, LV};
  \path (2.5,0) node[align=center](=B) {$\leftrightarrow$};
  \path (5,0) node[align=center,draw,rounded corners](FLV) {Adjoint Discretization\\ Forward PDE, LV};
  \path (5,-1) node[align=center](=R) {$\updownarrow$ $(\star )$};
  \path (5,-2) node[align=center,draw,rounded corners](FSLV) {Adjoint Discretization\\ Forward PDE, SLV};
  \path (2.5,-2) node[align=center](=U) {$\leftrightarrow$};
  \path (0,-2) node[align=center,draw,rounded corners](BSLV) {Discretization\\ Backward PDE, SLV};
\end{tikzpicture}
\end{center}
\vskip0.5cm
Here relationship $(\star )$ can only be achieved if similar discretizations are used for the backward PDEs stemming from the LV and SLV models.

An outline of the rest of our paper is as follows. 

In Section \ref{ForwardBackward} a relationship between the forward PDE and backward PDE is introduced, both for the case of the SLV model as for the case of the LV model.

In Section \ref{ForwardBackwardDiscrete} this relationship is preserved at the semidiscrete level: given a spatial discretization of the backward PDE, an adjoint spatial discretization for the forward PDE is defined such that both discretizations yield identical approximations for the fair value of non-path-dependent European options.

In Section \ref{BackwardForward} an actual spatial discretization, using second-order central finite difference schemes, is constructed for the backward PDE stemming from the SLV model and subsequently the corresponding adjoint spatial discretization is stated.

The main result of the paper is derived in Section \ref{LocalVolMatching}. 
It is shown that, under some assumptions, the adjoint spatial discretization can be employed to obtain an expression for the leverage function such that the approximation of the fair value of vanilla options is the same for the LV and SLV models.
In order to effectively use this expression, one has to solve a large system of non-linear ordinary differential equations (ODEs).

In Section \ref{TimeDiscretization} an Alternating Direction Implicit (ADI) temporal discretization scheme is applied to increase the computational efficiency in the numerical solution of this ODE system and in Section \ref{Calibration} an iteration procedure is described for handling the non-linearity.

In Section \ref{Numerical experiments} ample numerical experiments are presented to illustrate the performance of the obtained SLV calibration procedure.

The final Section \ref{Conclusion} gives concluding remarks.

\setcounter{equation}{0}
\section{Relationship between the forward and the backward Kolmogorov equation}\label{ForwardBackward}

Consider a European-style option with maturity $T$ and payoff $u_{0}$. 
Denote by $u(x,v,t)$ the \textit{non-discounted fair value} of the option under the SLV model \eqref{eq:SLVmodelX} at time to maturity $t$, that is at time level $\tau=T-t$, if $S_{\tau}=S_{0}\exp(x)$ and $V_{\tau} = v$.
It is well-known, see e.g.\ \cite{C11}, that the function $u$ satisfies the \textit{backward Kolmogorov equation}
\begin{equation}
\begin{array}{lll}
\tfrac{\partial}{\partial t} u &=& \tfrac{1}{2} \sigma^{2}_{SLV}(x,T-t)\psi^{2}(v) \tfrac{\partial^{2}}{\partial x^{2}}u  + \rho \xi \sigma_{SLV}(x,T-t)\psi(v) v^{\alpha} \tfrac{\partial^{2}}{\partial x \partial v} u + \tfrac{1}{2} \xi^{2} v^{2\alpha} \tfrac{\partial^{2}}{\partial v^{2}} u \\\\
&& + \ (r_{d}-r_{f}- \tfrac{1}{2} \sigma^{2}_{SLV}(x,T-t)\psi^{2}(v)) \tfrac{\partial}{\partial x} u + \kappa(\eta - v) \tfrac{\partial}{\partial v} u,
\end{array}
\label{eq:BackwardKolmogorov}
\end{equation}
for $x, v \in \R$, $0 < t \leq T$. At maturity, i.e.\ at time level $\tau = T$, the initial condition $u(x,v,0)$ is defined by the payoff $u_{0}$ of the option.
By solving PDE \eqref{eq:BackwardKolmogorov}, the fair value $e^{-r_{d}T}u(X_{0},V_{0},T)$ of the option under the SLV model can be determined at the spot, i.e.\ at $\tau =0$. For strictly positive values of the parameter $\alpha$, the process $V_{\tau}$ is non-negative and the spatial domain in the $v$-direction reduces to $v \geq 0$. 

If the option under consideration is non-path-dependent, then the payoff $u_{0}$ is only a function of $(X_{T},V_{T})$, the initial condition is given by $u(x,v,0) = u_{0}(x,v)$ and the non-discounted fair value $u(x,v,t)$ of the option can be written as
\begin{equation*}
u(x,v,t) = \mathbb{E}[u_{0}(X_{T},V_{T}) \vert X_{T-t} = x, V_{T-t} = v],
\end{equation*}
for $0 \le t \le T$.
By making use of the tower property for conditional expectations it readily follows that
\begin{equation}
u(X_{0},V_{0},T) = \mathbb{E}[ u(X_{T-t},V_{T-t},t) \vert X_{0}, V_{0} ] = \int_{-\infty}^{\infty} \int_{-\infty}^{\infty} u(x,v,t) p(x,v,T-t;X_{0},V_{0}) dx dv,
\label{eq:FairValue}
\end{equation}
for $0 \le t \le T$.
Recall that $p(x,v,\tau;X_{0},V_{0})$ denotes the joint density of $(X_{\tau},V_{\tau})$ under the SLV model \eqref{eq:SLVmodelX}. If the parameter $\alpha$ is chosen strictly positive, then the integral with respect to $v$ can be taken from $v=0$.
In particular, the fair value of non-path-dependent European options at the spot can also be computed by evaluating the integral
\begin{equation}
e^{-r_{d}T} \int_{-\infty}^{\infty} \int_{-\infty}^{\infty} u(x,v,0) p(x,v,T;X_{0},V_{0}) dx dv,
\label{eq:ForwardPricing}
\end{equation}
where $u(x,v,0)$ is defined by the payoff of the option. 

It can be shown, see e.g.\ \cite{R89}, that the joint density $p(x,v,\tau;X_{0},V_{0})$ satisfies the \textit{forward Kolmogorov equation}
\begin{equation}
\begin{array}{lll}
\tfrac{\partial}{\partial \tau} p &=& \tfrac{1}{2} \tfrac{\partial^{2}}{\partial x^{2}} \left( \sigma^{2}_{SLV}\psi^{2}(v)p \right) + \tfrac{\partial^{2}}{\partial x \partial v} \left( \rho \xi \sigma_{SLV}\psi(v)v^{\alpha} p \right) + \tfrac{1}{2} \tfrac{\partial^{2}}{\partial v^{2}} \left( \xi^{2} v^{2\alpha} p \right) \\\\
&& - \ \tfrac{\partial}{\partial x} \left( (r_{d}-r_{f}-\tfrac{1}{2}\sigma^{2}_{SLV}\psi^{2}(v)) p \right) - \tfrac{\partial}{\partial v} \left( \kappa(\eta - v) p \right),
\end{array}
\label{eq:ForwardKolmogorov}
\end{equation}
for $x, v \in \mathbb{R}, \tau >0$ and with initial condition $p(x,v,0;X_{0},V_{0}) = \delta(x-X_{0})\delta(v-V_{0})$ where $\delta$ denotes the Dirac delta function. For ease of presentation, the dependency of $\sigma_{SLV}$ on $(x,\tau)$ and the dependency of $p$ on $(x,v,\tau;X_{0},V_{0})$ is omitted in \eqref{eq:ForwardKolmogorov}.
Recall that the process $V_{\tau}$ is non-negative whenever $\alpha$ is strictly positive. In this case the spatial domain of the PDE in the $v$-direction is naturally restricted to $v \geq 0$. 
The integrals in \eqref{eq:CondExpec}, \eqref{eq:FairValue}, \eqref{eq:ForwardPricing} with respect to the $v$-variable can then be taken from $0$ to infinity. 

Equation \eqref{eq:FairValue} establishes a \textit{fundamental relationship between the forward and backward Kolmogorov equation}. 
It states that the fair value of non-path-dependent European options under the SLV model can be seen as the combination of the solution of two different PDEs. By considering the extreme time value $\tau = 0$ ($t=T$), or $\tau = T$ ($t=0$), only one PDE has to be solved.
In the forthcoming sections relationship \eqref{eq:FairValue} will be employed to define an adjoint spatial discretization for the forward equation.

Even if the functions $u$ and $p$ are known exactly, the integrals in \eqref{eq:FairValue} can often not be calculated analytically and one relies on numerical integration methods in order to approximate them.
In this article we assume that the integrand is known on a Cartesian grid.
Denote by $m_{1}$, respectively $m_{2}$, the number of spatial grid points in the $x$-direction, respectively $v$-direction.
The Cartesian grid is given by
\begin{equation}
(x_{i},v_{j}) \quad \mathrm{for} \ 1 \leq i \leq m_{1}, \ 1 \leq j \leq m_{2},
\label{eq:CartesianGrid}
\end{equation}
with $X_{\min}=x_{1} < x_{2} < \cdots < x_{m_{1}} = X_{\max}$, $V_{\min}=v_{1} < v_{2} < \cdots < v_{m_{2}} = V_{\max}$ and $X_{\min} < X_{0} < X_{\max}$, $V_{\min} < V_{0} < V_{\max}$. Define spatial mesh widths $\Delta x_{i} = x_{i}-x_{i-1}$ for $2\leq i \leq m_{1}$, $\Delta v_{j} = v_{j}-v_{j-1}$ for $2 \leq j \leq m_{2}$ and put $\Delta x_{1} = \Delta x_{m_{1}+1} = \Delta v_{1} = \Delta v_{m_{2}+1} = 0$.
When working with Cartesian grids, most numerical integration methods approximate the expression \eqref{eq:FairValue} by
\begin{equation}
u(X_{0},V_{0},T) \approx \sum_{i=1}^{m_{1}}\sum_{j=1}^{m_{2}} p(x_{i},v_{j},T-t;X_{0},V_{0}) u(x_{i},v_{j},t)  w_{x,i}w_{v,j},
\label{eq:DiscreteIntegraal}
\end{equation}
for certain weights $w_{x,i}, w_{v,j}$. If the numerical integration is performed with the trapezoidal rule, then the weights are given by
\begin{equation*}
w_{x,i} = \tfrac{\Delta x_{i} + \Delta x_{i+1}}{2} \quad \mathrm{for} \ 1 \le i \le m_{1}, \qquad w_{v,j} = \tfrac{\Delta v_{j} + \Delta v_{j+1}}{2} \quad \mathrm{for} \ 1 \le j \le m_{2}.
\end{equation*}
The values $X_{\min}, V_{\min}$, $X_{\max}, V_{\max}$ have to lie sufficiently far away from $(X_{0},V_{0})$ such that the truncation error is negligible.
If $\alpha$ is strictly positive, the value $V_{\min}$ can be set equal to zero. 
In order for \eqref{eq:DiscreteIntegraal} to be exact if $t=T$ $(\tau = 0)$, it is assumed that there exist indices $i_{0}, j_{0}$ such that $(x_{i_{0}},v_{j_{0}}) = (X_{0},V_{0})$ and the approximation $p(x,v,0;X_{0};V_{0}) \approx p_{0}(x,v)$ is used with
\begin{equation*}
p_{0}(x,v) = \left\{ \begin{array}{ll}
\tfrac{1}{w_{x,i_0} w_{v,j_0}} & \quad \mathrm{if} \ (x,v) \in [x_{i_{0}} - \tfrac{\Delta x_{i_{0}}}{2}, x_{i_{0}} + \tfrac{\Delta x_{i_{0}+1}}{2}] \times [v_{j_{0}} - \tfrac{\Delta v_{j_{0}}}{2}, v_{j_{0}} + \tfrac{\Delta v_{j_{0}+1}}{2}], \\\\
0 & \quad \mathrm{else.}
\end{array} \right.
\end{equation*}

Analogously as above, consider within the LV model a European-style option with payoff $u_{LV,0}$ at maturity $T$ and denote by $u_{LV}(x,t)$ the non-discounted fair value of the option under the LV model \eqref{eq:LVmodelX} at time $\tau = T-t$ if $X_{LV,\tau}=x$. The function $u_{LV}$ satisfies the backward Kolmogorov equation
\begin{equation}
\tfrac{\partial}{\partial t} u_{LV} = \tfrac{1}{2} \sigma^{2}_{LV} \tfrac{\partial^{2}}{\partial x^{2}} u_{LV} + (r_{d}-r_{f}-\tfrac{1}{2} \sigma^{2}_{LV}) \tfrac{\partial}{\partial x} u_{LV},
\label{eq:BackwardKolmogorovLV}
\end{equation}
for $x \in \R$, $0<t\leq T$ and with initial condition $u_{LV}(x,0)$ defined by the payoff $u_{LV,0}$ of the option. The non-discounted fair value of the option at the spot $(\tau = 0)$ is then given by $u_{LV}(X_{0},T)$.
For non-path-dependent options this fair value can also be formulated as
\begin{equation}
u(X_{0},T) = \int_{-\infty}^{\infty} u_{LV}(x,t) p_{LV}(x,T-t;X_{0}) dx,
\label{eq:FairValueLV}
\end{equation}
for $0 \le t \le T$, where $p_{LV}(x,\tau;X_{0})$ denotes the density of the process $X_{LV,\tau}$ in the LV model \eqref{eq:LVmodelX}.
It can be shown, see e.g.\ \cite{C11}, that this density function satisfies the forward Kolmogorov equation
\begin{equation} \tfrac{\partial}{\partial \tau} p_{LV} = \tfrac{1}{2} \tfrac{\partial^{2}}{\partial x^{2}} \left( \sigma^{2}_{LV}p_{LV} \right) -  \tfrac{\partial}{\partial x} \left( (r_{d}-r_{f}-\tfrac{1}{2}\sigma^{2}_{LV}) p_{LV} \right), 
\label{eq:ForwardKolmogorovLV}
\end{equation}
for $x \in \R, \tau >0$, and with initial condition $p_{LV}(x,0;X_{0})= \delta(x-X_{0})$.
Hence, the expression \eqref{eq:FairValueLV} establishes a fundamental relationship between the forward and backward Kolmogorov equation which is similar to \eqref{eq:FairValue}.
By applying the same numerical integration technique as above, the fair value from \eqref{eq:FairValueLV} can be approximated by
\begin{equation*}
u(X_{0},T) \approx \sum_{i=1}^{m_{1}}p_{LV}(x_{i},T-t;X_{0}) u(x_{i},t)  w_{x,i}.
\end{equation*}

Recall that the SLV model is calibrated perfectly to the LV model if the leverage function is defined by \eqref{eq:SigmatoMatch}.
It was shown by Gy\"ongy \cite{G86} that under this assumption both processes $X_{\tau}$ and $X_{LV,\tau}$ have the same marginal densities, i.e.\ that 
\begin{equation}
\int_{-\infty}^{\infty} p(x,v,\tau;X_{0},V_{0})dv = p_{LV}(x,\tau;X_{0})
\label{eq:SamePriceContinuous}
\end{equation} 
for $x \in \R, \tau \geq 0.$
From now on, for the ease of presentation, the dependency of $p$ and $p_{LV}$ on the spot values $X_{0}, V_{0}$ is omitted.

\setcounter{equation}{0}
\section{Adjoint spatial discretization}\label{ForwardBackwardDiscrete}

In general the values $p(x_{i},v_{j},\tau)$ and $u(x_{i},v_{j},t)$ are not known exactly if $\tau>0$ and $t>0$, respectively, and one relies on numerical methods to approximate them. 
Extensive literature is available on numerical techniques to solve backward Kolmogorov equations, see e.g.\ \cite{TR00}.
A common approach in financial mathematics in order to approximate the fair value of options is by numerically solving the pertinent PDE  using the general \textit{method of lines} (MOL), cf.\ \cite{HV03}. 
In this approach, the PDE is first discretized in the spatial variables $x$ and $v$, yielding a large system of stiff ordinary differential equations (ODEs). This, so-called, \textit{semidiscrete system} is subsequently solved by applying a suitable implicit time stepping method. 

Spatial discretization by finite difference or finite volume methods of the backward Kolmogorov equation \eqref{eq:BackwardKolmogorov} on a Cartesian grid \eqref{eq:CartesianGrid} yields approximations $\boldsymbol{U}_{i,j}(t)$ of the exact non-discounted option value $u(x_{i},v_{j},t)$. 
Denote by $\boldsymbol{U}(t)$ the $m_{1} \times m_{2}$ matrix with entries $\boldsymbol{U}_{i,j}(t)$ and denote by $\boldsymbol{P}(\tau)$ a matrix with entries $\boldsymbol{P}_{i,j}(\tau)$ that represent approximations to the exact density values $p(x_{i},v_{j},\tau)$.
In this section, for a general spatial discretization of the backward Kolmogorov equation, an adjoint spatial discretization of the corresponding forward equation is defined such that 
\begin{equation}
\sum_{i=1}^{m_{1}}\sum_{j=1}^{m_{2}} \boldsymbol{P}_{i,j}(T-t) \boldsymbol{U}_{i,j}(t) w_{x,i}w_{v,j}
\label{eq:DiscreteSom}
\end{equation}
is constant for $0 \le t \le T$. This can be viewed as a discrete version of relationship \eqref{eq:FairValue}.

Let the vector
$$ U(t) = \mathrm{vec}[\boldsymbol{U}(t)],$$
where $\mathrm{vec}[\cdot]$ denotes the operator that turns any given matrix into a vector by putting its successive columns below each other. 
Spatial discretization of \eqref{eq:BackwardKolmogorov} leads to a large system of ODEs of the form
\begin{equation}
U'(t) = A^{(B)}(t)U(t),
\label{eq:StandardBackwardDiscretization}
\end{equation} 
for $0 < t \leq T$, with given matrix $A^{(B)}(t)$ and with given vector $U(0)$ that is defined by the payoff of the option.
Let the vector 
$$P(\tau) = \mathrm{vec}[\boldsymbol{P}(\tau)],$$
where the matrix $\boldsymbol{P}(0)$ is defined by the function $p_{0}$ from Section \ref{ForwardBackward}.
Note that it has only one non-zero entry.
Denote by $M$ the diagonal matrix with diagonal entries
$$ M_{k,k} = w_{x,i}w_{v,j}, $$
where $i,j$ are the indices such that element $U_{k}(t)$, respectively $P_{k}(\tau)$, corresponds to $\boldsymbol{U}_{i,j}(t)$, respectively $\boldsymbol{P}_{i,j}(\tau)$.
The semidiscrete analogue \eqref{eq:DiscreteSom} of \eqref{eq:DiscreteIntegraal} can then be compactly written as
\begin{equation}
u(X_{0},V_{0},T) \approx P(T-t)^{\rm{T}}MU(t),
\label{eq:SemidiscreteOptionPrice}
\end{equation}
where $^{\rm{T}}$ denotes taking the transpose.
If $k_{0}$ is the index that corresponds to $(i_{0},j_{0})$, then for $t=T$ the right-hand side of \eqref{eq:SemidiscreteOptionPrice} is equal to $U_{k_{0}}(T) = \boldsymbol{U}_{i_{0},j_{0}}(T)$. Now, consider the fair value of any given non-path-dependent European option with maturity $T$. It is readily seen that semidiscretization of the forward equation \eqref{eq:ForwardKolmogorov} and semidiscretization of the backward equation \eqref{eq:BackwardKolmogorov} define {\it the same approximation} \eqref{eq:SemidiscreteOptionPrice} of the fair value for all $0 \leq t \leq T$, i.e.\ property \eqref{eq:FairValue} holds in the semidiscrete sense, if 
\begin{equation}
U_{k_{0}}(T) = P(0)^{\rm{T}}MU(T) = P(T-t)^{\rm{T}}MU(t) = P(\tau)^{\rm{T}}MU(T-\tau)
\label{eq:SemiForwardMatchBackward}
\end{equation}
for all $0\leq t, \tau \leq T$. This requirement is satisfied whenever 
\begin{equation*}
0 = P'(\tau)^{\rm{T}} M U(T-\tau) -  P(\tau)^{\rm{T}}MA^{(B)}(T-\tau)U(T-\tau),
\end{equation*}
holds for all $0 \leq \tau \leq T$.
Accordingly, we define the \textit{adjoint spatial discretization} of the forward Kolmogorov equation \eqref{eq:ForwardKolmogorov} as 
\begin{equation}
P'(\tau) = M^{-1}(A^{(B)}(T-\tau))^{\rm{T}}MP(\tau) \quad \mathrm{for} \ 0 \leq \tau \leq T,
\label{eq:DefinitieVoorwDiscr}
\end{equation}
In this paper we always employ the adjoint spatial discretization \eqref{eq:DefinitieVoorwDiscr} of the forward equation. Thus, given any semidisretization of the backward equation, the obtained approximated option value for non-path-dependent European options satisfy \eqref{eq:SemiForwardMatchBackward}.

It is convenient to introduce
\begin{equation}
\overline{P}(\tau) = MP(\tau).
\label{eq:MP}
\end{equation}
Denote by $\overline{\boldsymbol{P}}(\tau)$ the matrix corresponding to the vector $\overline{P}(\tau)$. 
The elements $\overline{\boldsymbol{P}}_{i,j}(\tau)$ can be viewed as approximations of
$$ \int_{x_{i} - \tfrac{\Delta x_{i}}{2}}^{x_{i} + \tfrac{\Delta x_{i+1}}{2}} \int_{v_{j} - \tfrac{\Delta v_{j}}{2}}^{v_{j} + \tfrac{\Delta v_{j+1}}{2}} p(x,v,\tau)dx dv, $$
and hence, as an approximation of the probability that
$$ (X_{\tau},V_{\tau}) \in [x_{i} - \tfrac{\Delta x_{i}}{2}, x_{i} + \tfrac{\Delta x_{i+1}}{2}] \times [v_{j} - \tfrac{\Delta v_{j}}{2}, v_{j} + \tfrac{\Delta v_{j+1}}{2}]. $$
Clearly,
\begin{equation*}
\overline{P}'(\tau) = MM^{-1} (A^{(B)}(T-\tau))^{\rm{T}}M P(\tau) = (A^{(B)}(T-\tau))^{\rm{T}}\overline{P}(\tau),
\end{equation*}
for $0 \leq \tau \leq T$, with given vector $\overline{P}(0) = MP(0)$, i.e.\ with
\begin{equation*}
\overline{P}_{k}(0) = \left\{ \begin{array}{ll}
1 & \qquad \mathrm{if} \ k = k_{0}, \\
0 & \qquad \mathrm{else.}
\end{array} \right.
\end{equation*}
The approximation \eqref{eq:SemidiscreteOptionPrice} of the non-discounted fair value of a non-path-dependent European option at the spot can then also be represented as
$$ u(X_{0},V_{0},T) \approx U_{k_{0}}(T) = \overline{P}(T-t)^{\rm{T}}U(t).$$

\setcounter{equation}{0}
\section{Spatial discretization by finite differences}\label{BackwardForward}

In this section, a spatial discretization of the backward equation \eqref{eq:BackwardKolmogorov} by finite differences will be performed on a non-uniform Cartesian grid \eqref{eq:CartesianGrid}. This semidiscretization then defines the adjoint spatial discretization of the forward equation \eqref{eq:ForwardKolmogorov}.

\subsection{Spatial discretization of the backward equation}
\label{subsec:BackwardDiscretization}

To construct a spatial grid and a semidiscretization for \eqref{eq:BackwardKolmogorov}, the spatial domain needs to be truncated to a bounded set $[X_{\min}, X_{\max}] \times [V_{\min}, V_{\max}]$. The boundaries have to lie sufficiently far away from $(X_{0},V_{0})$ such that the truncation error incurred is negligible. 
If the parameter $\alpha$ is strictly positive, the process $V_{\tau}$ is non-negative and $V_{\min}$ is naturally set equal to zero.
For non-path-dependent European options the following boundary conditions are imposed:
\begin{equation}
\begin{array}{ll}
\tfrac{\partial^{2}}{\partial x^{2}} u(X_{\min},v,t) = \tfrac{\partial}{\partial x} u(X_{\min},v,t) & \qquad \mathrm{for} \quad 0 \leq v \leq V_{\max}, \ 0 < t \leq T, \\\\
\tfrac{\partial^{2}}{\partial x^{2}} u(X_{\max},v,t) = \tfrac{\partial}{\partial x} u(X_{\max},v,t) & \qquad \mathrm{for} \quad 0 \leq v \leq V_{\max}, \ 0 < t \leq T.
\end{array}
\label{eq:BCBackwardX}
\end{equation}
The above conditions at $x=X_{\min}$ and $x=X_{\max}$ correspond to linear boundary conditions in the $s$-variable, where $s=S_{0}\exp(x)$, cf.\ Section \ref{intro}.
If $\alpha = 0$, the process $V_{\tau}$ can take negative values and it is additionally assumed that
\begin{equation*}
\begin{array}{ll}
\tfrac{\partial^{2}}{\partial v^{2}} u(x,V_{\min},t) = 0 & \qquad \mathrm{for} \quad X_{\min} \leq x \leq X_{\max}, \ 0 < t \leq T, \\\\
\tfrac{\partial^{2}}{\partial v^{2}} u(x,V_{\max},t) = \tfrac{\partial}{\partial v} u(x,V_{\max},t) & \qquad \mathrm{for} \quad X_{\min} \leq x \leq X_{\max}, \ 0 < t \leq T.
\end{array}
\end{equation*}
Thus at $v=V_{\min}$ a linear boundary condition is taken. The condition at $v=V_{\max}$ corresponds with a linear boundary condition in the variable $\exp(v)$.
For values of $\alpha$ that are strictly positive, the process $V_{\tau}$ is non-negative. Moreover, $V_{\tau}=0$ can be attained if $0< \alpha \le 1/2$. In these cases the boundary $v=0$ of the PDE requires special attention.
It has recently been proved in \cite{ET11} that setting $v=0$ in the PDE \eqref{eq:BackwardKolmogorov} at this boundary then yields the correct condition here. At the boundary $v=V_{\max}$ it is then additionally assumed that
\begin{equation*}
\tfrac{\partial^{2}}{\partial v^{2}} u(x,V_{\max},t) = 0 \qquad \mathrm{for} \quad X_{\min} \leq x \leq X_{\max}, \ 0 < t \leq T.
\end{equation*} 

\begin{figure}
\begin{center}
\includegraphics[scale=0.5]{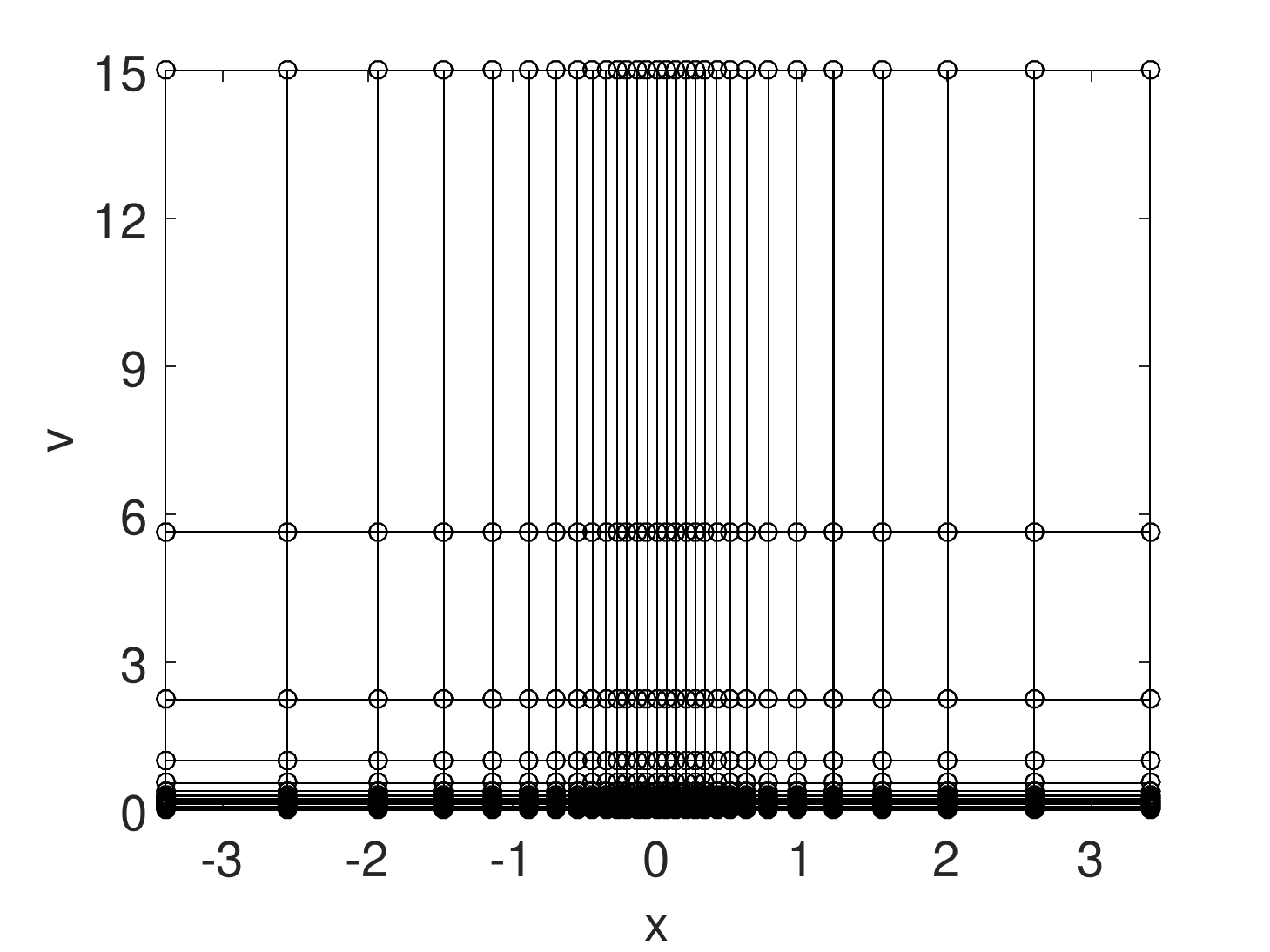} 
\includegraphics[scale=0.5]{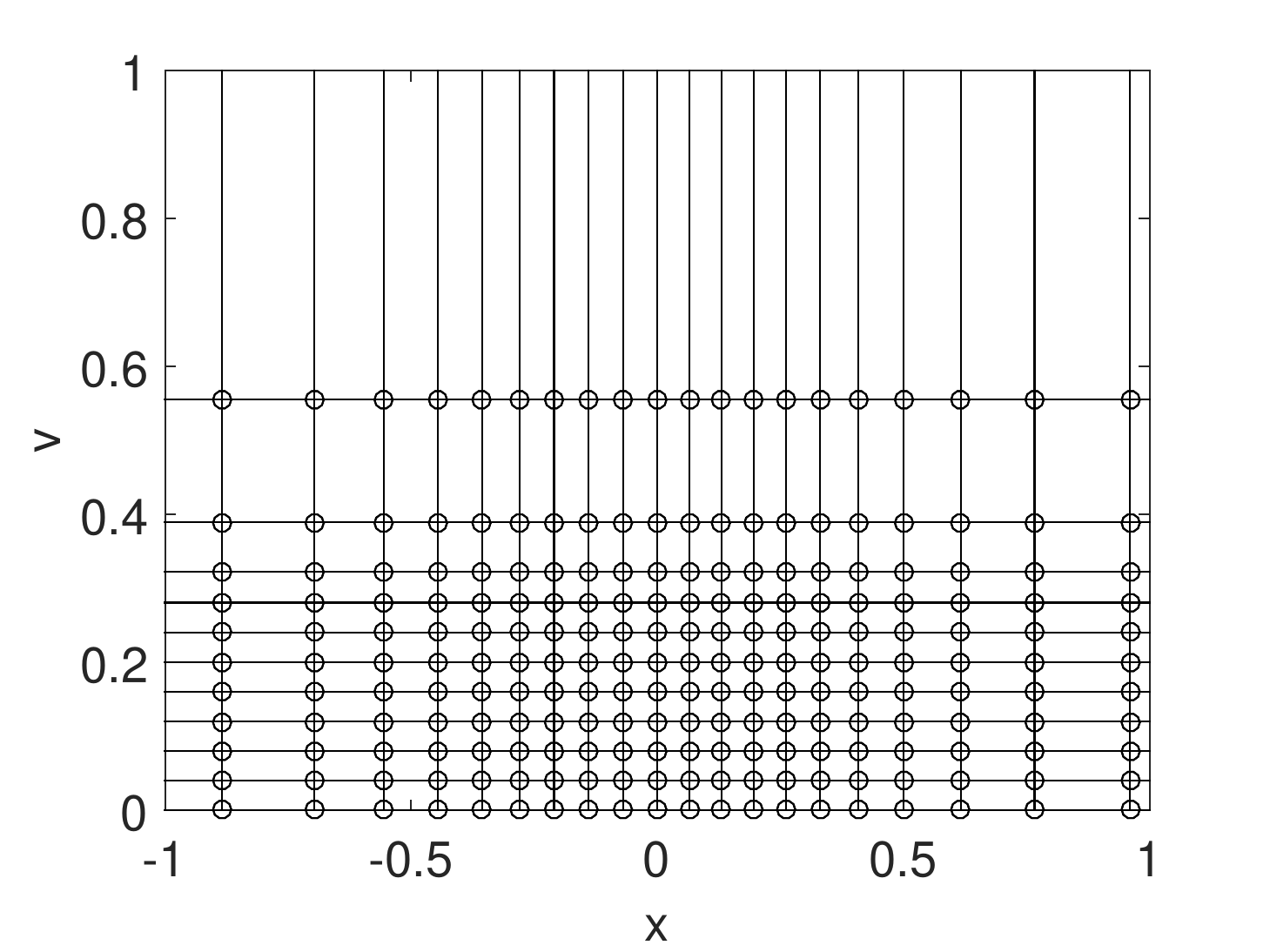}
\caption{Sample grid for $m_{1}=30, m_{2}=15$, the case $\alpha >0$, $X_{\min} = - \log(30), X_{\max} = \log(30)$, $V_{\min}=0, V_{\max}=15$ and $(X_{0},V_{0})=(0,0.2)$. The left plot displays the complete grid, the right plot shows the uniformity around $(X_{0},V_{0})$.}
\label{fig:SpatialGrid}
\end{center}
\end{figure}

For this truncated domain, non-uniform meshes are applied in both the $x$- and $v$-direction such that relatively many mesh points lie in the neighbourhood of $x=X_{0}$ and $v=V_{0}$. The application of such non-uniform meshes improves the accuracy of the finite difference discretization compared to using uniform meshes. 
The type of non-uniform meshes that is employed is similar to the ones considered in e.g.\ \cite{HH12,IHF10,TR00}, which make use of a transformation with the $\sinh$ function of an underlying uniform mesh. 
Akin to \cite{HH12} we opt to have here a fine, locally uniform mesh around the point $(x,v) = (X_{0},V_{0})$. 
Further, the choice $m_{1}=2m_{2}$ is considered. 
Denote by $\Delta \xi$ the mesh width of the underlying uniform mesh in the $x$-direction. 
Then the mesh under consideration in the $x$-direction is \textit{smooth} in the sense that there exist real constants $C_{0}, C_{1}, C_{2} > 0$ such that the mesh widths $\Delta x_{i}=x_i-x_{i-1}$ satisfy
\begin{equation*}
C_{0} \Delta \xi \le \Delta x_{i} \le C_{1} \Delta \xi \qquad \mathrm{and} \qquad \vert \Delta x_{i+1} - \Delta x_{i} \vert \le C_{2} (\Delta \xi)^{2}
\end{equation*} 
uniformly in $i$ and $m_{1}$. Analogously, the mesh in the $v$-direction is smooth.
As an illustration, the left plot in Figure \ref{fig:SpatialGrid} displays the spatial grid for the (small) sample values $m_{1}=30, m_{2}=15$, in the case $\alpha>0$, $X_{\min} = - \log(30), X_{\max} = \log(30)$, $V_{\min}=0, V_{\max}=15$ and $(X_{0},V_{0}) = (0,0.2)$. The right plot in Figure \ref{fig:SpatialGrid} displays a part of the spatial grid to show the local uniformity of the grid around $(X_{0},V_{0})$.

Semidiscretization on the spatial grid is performed by finite differences. Let $f: \R \rightarrow \R$ be any given function.
To approximate the first derivative $f'(x_{i})$ we consider two finite difference schemes:
\begin{subeqnarray}
f'(x_{i}) &\approx & \beta_{i,-1} f(x_{i-1}) + \beta_{i,0} f(x_{i}) + \beta_{i,1} f(x_{i+1}), \slabel{eq:FirstDerivCentral} \\
f'(x_{i}) &\approx & \gamma_{i,0} f(x_{i}) + \gamma_{i,1} f(x_{i+1}) + \gamma_{i,2} f(x_{i+2}), \slabel{eq:FirstDerivForw}
\end{subeqnarray} 
with coefficients given by
\begin{equation*}
\begin{array}{lll}
\beta_{i,-1} = \tfrac{-\Delta x_{i+1}}{\Delta x_{i} (\Delta x_{i} + \Delta x_{i+1})}, & \beta_{i,0} = \tfrac{\Delta x_{i+1} - \Delta x_{i}}{\Delta x_{i} \Delta x_{i+1}}, & \beta_{i,1} = \tfrac{\Delta x_{i}}{\Delta x_{i+1} (\Delta x_{i} + \Delta x_{i+1})}, \\\\
\gamma_{i,0} = \tfrac{-2\Delta x_{i+1} - \Delta x_{i+2}}{\Delta x_{i+1} (\Delta x_{i+1} + \Delta x_{i+2})}, & \gamma_{i,1} = \tfrac{\Delta x_{i+1} + \Delta x_{i+2}}{\Delta x_{i+1} \Delta x_{i+2}}, & \gamma_{i,2} = \tfrac{-\Delta x_{i+1}}{\Delta x_{i+2} (\Delta x_{i+1} + \Delta x_{i+2})}. 
\end{array}
\end{equation*}
To approximate the second derivative $f''(x_{i})$, we apply the finite difference scheme
\begin{equation}
f''(x_{i}) \approx \delta_{i,-1} f(x_{i-1}) + \delta_{i,0} f(x_{i}) + \delta_{i,1} f(x_{i+1}),
\label{eq:SecondDerivCentral}
\end{equation}
where
$$ \delta_{i,-1} = \tfrac{2}{\Delta x_{i} (\Delta x_{i} + \Delta x_{i+1})}, \quad \delta_{i,0} = \tfrac{-2}{\Delta x_{i} \Delta x_{i+1}}, \quad \delta_{i,1} = \tfrac{2}{\Delta x_{i+1} (\Delta x_{i} + \Delta x_{i+1})}. $$
The finite difference schemes in the $v$-direction are defined completely analogously.
For a function of two variables $f: \R^{2} \rightarrow \R $, the mixed derivative is approximated by application of \eqref{eq:FirstDerivCentral}~successively in the two directions.
The finite difference schemes above are all well-known. Formula \eqref{eq:FirstDerivCentral}, respectively \eqref{eq:FirstDerivForw}, is called the second-order central, respectively second-order forward, formula for the first derivative. Finite difference scheme \eqref{eq:SecondDerivCentral} is called the second-order central formula for the second derivative.
Through Taylor expansion it can be verified that each of the finite difference approximations above has a second-order truncation error, provided that the function $f$ is sufficiently often continuously differentiable and the mesh is smooth, cf.\ above.

The actual semidiscretization of the backward PDE \eqref{eq:BackwardKolmogorov} is defined as follows.
At all spatial grid points that do not lie on the boundary of the truncated domain, each spatial derivative appearing in \eqref{eq:BackwardKolmogorov} is replaced by its corresponding second-order central finite difference scheme described above. 

Concerning the boundaries in the $x$-direction, it is assumed that the pertinent conditions from \eqref{eq:BCBackwardX} are valid for every $x$ smaller than $x_{2}$ or larger than $x_{m_{1}-1}$. Thus for these extreme values $x$ we assume that $u(\cdot,v,t)$ is an exponential function. For instance, considering the upper boundary
$$ u_{x} = u_{xx} = C_{b,1} \exp(x) \qquad \mathrm{whenever} \quad x > x_{m_{1}-1},$$
and hence
$$ u = C_{b,1} \exp(x) + C_{b,2} \qquad \mathrm{whenever} \quad x > x_{m_{1}-1},$$
for some constants $C_{b,1}$, $C_{b,2}$.
Based on the function values $u$ at $x_{m_{1}-1}$ and $x_{m_{1}}$ the constants $C_{b,1}$, $C_{b,2}$ can be determined and this leads to the approximation for both the first and second derivatives in the $x$-direction at $x_{m_{1}}$ given by
$$ 
C_{b,1} \exp(x_{m_{1}}) =
\tfrac{-\exp(x_{m_{1}})}{\exp(x_{m_{1}})-\exp(x_{m_{1}-1})} u(x_{m_{1}-1},v,t) +  \tfrac{\exp(x_{m_{1}})}{\exp(x_{m_{1}})-\exp(x_{m_{1}-1})} u(x_{m_{1}},v,t).
$$
The first and second derivatives in the $x$-direction at the lower boundary are approximated analogously.

If the parameter $\alpha$ is strictly positive, a linear boundary condition is applied at $v=V_{\max}$, i.e.\ the second derivative in the $v$-direction is equal to zero. The first derivative in the $v$-direction at this boundary is approximated by using the central scheme \eqref{eq:FirstDerivCentral} with the virtual point $V_{\max}+\Delta v_{m_{2}}$, where the value at this point is defined by extrapolation using the linear boundary condition. This discretization reduces to the first-order backward finite difference formula for the first derivative. Moreover, at the boundary $V_{\min}=0$ all the second derivatives vanish and the first derivative $\partial u / \partial v $ is then approximated by using the forward scheme \eqref{eq:FirstDerivForw}.
If $\alpha$ equals zero, a linear boundary condition is imposed at $V_{\min}$ and discretization of the spatial derivatives in the $v$-direction is performed as above. The discretization of the boundary condition at $V_{\max}$ is then performed analogously as the discretization of the boundary conditions in the $x$-direction.

Denote by $D_{x}$, respectively $D_{xx}$, the matrices corresponding with the first, respectively second, derivatives in the $x$-direction. Analogously, denote by $D_{v}, D_{vv}$ the matrices corresponding with spatial derivatives in the $v$-direction. For example, given that $\alpha >0$, $D_{v}$ is the matrix with entries
\begin{equation*}
\begin{array}{lll}
(D_{v})_{1,1} = \tfrac{-2\Delta v_{2} - \Delta v_{3}}{\Delta v_{2} (\Delta v_{2} + \Delta v_{3})}, & (D_{v})_{1,2} = \tfrac{\Delta v_{2} + \Delta v_{3}}{\Delta v_{2} \Delta v_{3}}, & (D_{v})_{1,3} = \tfrac{-\Delta v_{2}}{\Delta v_{3} (\Delta v_{2} + \Delta v_{3})}, \\\\
(D_{v})_{i,i-1} = \tfrac{-\Delta v_{i+1}}{\Delta v_{i} (\Delta v_{i} + \Delta v_{i+1})}, & (D_{v})_{i,i} = \tfrac{\Delta v_{i+1} - \Delta v_{i}}{\Delta v_{i} \Delta v_{i+1}}, & (D_{v})_{i,i+1} = \tfrac{\Delta v_{i}}{\Delta v_{i+1} (\Delta v_{i} + \Delta v_{i+1})}, \\\\
(D_{v})_{m_{2},m_{2}-1} = \tfrac{-1}{\Delta v_{m_{2}}}, & (D_{v})_{m_{2},m_{2}} = \tfrac{1}{\Delta v_{m_{2}}},&
\end{array}
\end{equation*}
where $i \in \{2,\ldots,m_{2}-1\}$. Denote by $L(\tau)$ the $m_{1} \times m_{1}$ diagonal matrix with entries $\sigma_{SLV}(x_{i},\tau)$, let $\Lambda$ be the $m_{2} \times m_{2}$ diagonal matrix with entries $v_{j}$, and define for an arbitrary function $f:\R \rightarrow \R$ the matrix $f(\Lambda)$ as the diagonal matrix with entries $f(v_{j})$. Further, denote by $I_{x}$, respectively $I_{v}$, the identity matrix of size $m_{1} \times m_{1}$, respectively $m_{2} \times m_{2}$.
Then semidiscretization of \eqref{eq:BackwardKolmogorov} yields a system of differential equations given by
\begin{eqnarray*}
\boldsymbol{U}'(t) &=& \tfrac{1}{2} L^{2}(T-t) D_{xx} \boldsymbol{U}(t) \psi^{2}(\Lambda) + \rho \xi L(T-t) D_{x} \boldsymbol{U}(t) D^{\rm{T}}_{v} \Lambda^{\alpha} \psi(\Lambda) + \tfrac{1}{2} \xi^{2} \boldsymbol{U}(t) D^{\rm{T}}_{vv} \Lambda^{2\alpha} \\
&& + \ (r_{d}-r_{f}) D_{x} \boldsymbol{U}(t) - \tfrac{1}{2} L^{2}(T-t) D_{x} \boldsymbol{U}(t) \psi^{2}(\Lambda) + \boldsymbol{U}(t) D^{\rm{T}}_{v} \kappa(\eta I_{v} - \Lambda),
\end{eqnarray*}
for $0 < t \leq T$.
This can be written in the form \eqref{eq:StandardBackwardDiscretization}, 
\begin{equation}
U'(t) = A^{(B)}(t)U(t) = (A^{(B)}_{0}(t) + A^{(B)}_{1}(t) + A^{(B)}_{2}(t))U(t)
\label{eq:SemidiscretizationBackwardEquation}
\end{equation} 
for $0 < t \leq T$ where, using a well-known property of the \textit{Kronecker product} $\otimes$, 
\begin{eqnarray*}
A^{(B)}_{0}(t) &=& (\rho \xi \psi(\Lambda) \Lambda^{\alpha} D_{v}) \otimes (L(T-t) D_{x}), \\
A^{(B)}_{1}(t) &=& \tfrac{1}{2} \psi^{2}(\Lambda) \otimes (L^{2}(T-t) D_{xx}) + (r_{d}-r_{f}) I_{v} \otimes D_{x} - \tfrac{1}{2} \psi^{2}(\Lambda) \otimes (L^{2}(T-t) D_{x}), \\
A^{(B)}_{2}(t)    &=& ( \tfrac{1}{2} \xi^{2} \Lambda^{2\alpha} D_{vv} + \kappa(\eta I_{v} - \Lambda)D_{v}) \otimes I_{x}.
\end{eqnarray*}
The initial vector $U(0)$ is defined by the payoff of the option.

\subsection{Spatial discretization of the forward equation}

As indicated in Section \ref{ForwardBackwardDiscrete}, semidiscretization of the forward equation \eqref{eq:ForwardKolmogorov} is performed by the adjoint spatial discretization \eqref{eq:DefinitieVoorwDiscr}.
Since the transpose of the Kronecker product of two matrices is equal to the Kronecker product of the transposed matrices, and recalling that $t = T-\tau$, it follows that $\overline{P}(\tau)$ defined by \eqref{eq:MP} is given by the system of ODEs
\begin{equation} 
\overline{P}'(\tau) = \overline{A}^{(F)}(\tau) \overline{P}(\tau) = (\overline{A}^{(F)}_{0}(\tau) + \overline{A}^{(F)}_{1}(\tau) + \overline{A}^{(F)}_{2}(\tau))\overline{P}(\tau),
\label{eq:SemidiscretizationForwardEquation} 
\end{equation}
for $\tau > 0$, with
\begin{eqnarray*}
\overline{A}^{(F)}_{0}(\tau) &=& (\rho \xi D^{\rm{T}}_{v} \Lambda^{\alpha}\psi(\Lambda)) \otimes (D^{\rm{T}}_{x} L(\tau)), \\
\overline{A}^{(F)}_{1}(\tau) &=& \tfrac{1}{2} \psi^{2}(\Lambda) \otimes (D^{\rm{T}}_{xx} L^{2}(\tau)) + (r_{d}-r_{f}) I_{v} \otimes D^{\rm{T}}_{x} - \tfrac{1}{2} \psi^{2}(\Lambda) \otimes (D^{\rm{T}}_{x} L^{2}(\tau)), \\
\overline{A}^{(F)}_{2}(\tau) &=& ( \tfrac{1}{2} \xi^{2} D^{\rm{T}}_{vv} \Lambda^{2\alpha} + D^{\rm{T}}_{v} \kappa(\eta I_{v} - \Lambda)) \otimes I_{x},
\end{eqnarray*}
and given initial vector $\overline{P}(0)$.
This in turn corresponds to the system of differential equations 
\begin{eqnarray}
\overline{\boldsymbol{P}}'(\tau) &=& \tfrac{1}{2} D^{\rm{T}}_{xx} L^{2}(\tau) \overline{\boldsymbol{P}}(\tau) \psi^{2}(\Lambda) + \rho \xi D^{\rm{T}}_{x} L(\tau) \overline{\boldsymbol{P}}(\tau) \psi(\Lambda)\Lambda^{\alpha} D_{v}  + \tfrac{1}{2} \xi^{2} \overline{\boldsymbol{P}}(\tau) \Lambda^{2\alpha} D_{vv} \nonumber\\
&& + \ (r_{d}-r_{f}) D^{\rm{T}}_{x} \overline{\boldsymbol{P}}(\tau) - \tfrac{1}{2} D^{\rm{T}}_{x} L^{2}(\tau) \overline{\boldsymbol{P}}(\tau) \psi^{2}(\Lambda) + \overline{\boldsymbol{P}}(\tau) \kappa(\eta I_{v} - \Lambda) D_{v},
\label{eq:Pmatrix}
\end{eqnarray}
for $\tau > 0$.
The expression \eqref{eq:Pmatrix} shall be employed to calibrate the SLV model to the LV model.

The total integral of a density function is equal to one. For a natural adjoint spatial discretization \eqref{eq:SemidiscretizationForwardEquation} one would expect that the total numerical integral of $\boldsymbol{P}$, corresponding with $\overline{\boldsymbol{P}}$, is close to one.
Let $e_x$ and $e_v$ denote the vectors consisting of all ones with lengths $m_1$ and $m_2$, respectively.
By construction of the finite difference discretization and the chosen boundary conditions for the SLV model \eqref{eq:BackwardKolmogorov} there holds
\begin{equation}
D_{xx} e_x = D_{x} e_x = 0 \quad {\rm and} \quad D_{vv} e_v =D_{v} e_v = 0
\label{eq:DiscretizationBCV}
\end{equation}
and it directly follows that
$$ 
e_x^{\rm T}\, \overline{\boldsymbol{P}}'(\tau) e_v = 0 
$$
for all $\tau > 0$.
Since further $e_x^{\rm T}\, \overline{\boldsymbol{P}}(0) e_v = 1$, this yields
\begin{equation*}
\sum_{i=1}^{m_{1}}\sum_{j=1}^{m_{2}} \boldsymbol{P}_{i,j}(\tau) w_{x,i}w_{v,j} = \sum_{i=1}^{m_{1}}\sum_{j=1}^{m_{2}} \overline{\boldsymbol{P}}_{i,j}(\tau) = 1
\end{equation*}
for all $\tau \ge 0$. It can be concluded that the adjoint spatial discretization of the forward Kolmogorov equation keeps the total numerical integral of the density identically equal to one, which is a favourable property.

\setcounter{equation}{0}
\section{Matching the semidiscrete LV and SLV models}\label{LocalVolMatching}

In this section the main result of the article is presented. It is shown that, under some assumptions, one can calibrate the semidiscrete SLV model exactly to the corresponding semidiscrete LV model.

A priori the leverage function $\sigma_{SLV}$, and hence the matrix function $L$, are unknown and one wishes to determine them in such a way that the LV model and the SLV model define identical values for European call and put options.
In practice, however, even the LV function $\sigma_{LV}$ is not known analytically in general and one relies on numerical methods to approximate the option value defined by the LV model.
Accordingly, it is unrealistic to require an algorithm to produce a leverage function $\sigma_{SLV}$ such that the SLV model yields the same exact European call and put values as the LV model. 
One rather wants to construct the leverage function in such a way that the two models yield {\it identical approximate values}\, for European call and put options whenever similar semidiscretizations of these models are used.

A common approach to approximate the fair value of a European-style option under the LV model is by discretizing the backward PDE \eqref{eq:BackwardKolmogorovLV} with finite differences.
Since the region of interest in the $x$-direction in the LV model is the same as that in the SLV model, the same spatial mesh can be used in this spatial direction.
Denote by $U_{LV}(t)$ the vector with approximations $U_{LV,i}(t)$ to $u_{LV}(x_{i},t)$ for $1 \leq i \leq m_{1}$ such that the component $U_{LV,i_{0}}(T)$ is the approximation of the non-discounted fair value at the spot.
Semidiscretization by finite differences then leads to a system of ODEs 
\begin{equation}
U'_{LV}(t) = A^{(B)}_{LV}(t) U_{LV}(t) 
\label{eq:SemidiscretizationBackwardEquationLV}
\end{equation}
for $0 < t \leq T$, with initial vector $U_{LV}(0)$ defined by the payoff $u_{LV,0}$.
Since the spatial derivatives $\partial/\partial x$ and $\partial^{2}/\partial x^{2}$ in \eqref{eq:BackwardKolmogorovLV} also occur in the backward equation \eqref{eq:BackwardKolmogorov}, and since also the same boundary conditions from Section \ref{BackwardForward} can be applied for non-path-dependent European options, the same finite difference matrices $D_{x}$ and $D_{xx}$ can be used to perform semidiscretization, and hence 
$$ 
A^{(B)}_{LV}(t) = \tfrac{1}{2} L^{2}_{LV}(T-t) D_{xx} + (r_{d}-r_{f}) D_{x} - \tfrac{1}{2} L^{2}_{LV}(T-t) D_{x}, $$
where $L_{LV}(\tau)$ is the $m_{1} \times m_{1}$ diagonal matrix with entries $\sigma_{LV}(x_{i},\tau)$.

Denote by $P_{LV}(\tau)$ a vector with approximations $P_{LV,i}(\tau)$ of $p_{LV}(x_{i},\tau)$, where $p_{LV}$ is given by the forward equation \eqref{eq:ForwardKolmogorovLV}, and define $\overline{P}_{LV}(\tau)= M_{LV}P_{LV}(\tau)$ where $M_{LV}$ is the diagonal matrix with entries 
$$ (M_{LV})_{i,i} = w_{x,i}, \qquad \mathrm{for} \ 1 \leq i \leq m_{1}. $$
Analogously to Section \ref{ForwardBackwardDiscrete}, we define an adjoint forward discretization by
\begin{equation}
\overline{P}'_{LV}(\tau) = (A^{(B)}_{LV}(T-\tau))^{\rm{T}}\,\overline{P}_{LV}(\tau) = \overline{A}^{(F)}_{LV}(\tau)\overline{P}_{LV}(\tau)
\label{eq:SemidiscretizationForwardEquationLV}
\end{equation} 
for $\tau > 0$, with
$$ \overline{A}^{(F)}_{LV}(\tau) = \tfrac{1}{2} D^{\rm{T}}_{xx} L^{2}_{LV}(\tau) + (r_{d}-r_{f})D^{\rm{T}}_{x} - \tfrac{1}{2} D^{\rm{T}}_{x} L^{2}_{LV}(\tau),$$
and
$$ \overline{P}_{LV,i}(0) = \left\{ \begin{array}{ll}
1 & \quad \mathrm{if} \ i = i_{0}, \\
0 & \quad \mathrm{else,}
\end{array} \right. $$
so that
$$ U_{LV,i_{0}}(T) = \overline{P}_{LV}(T-t)^{\rm{T}}U_{LV}(t)$$
for all $0\leq t \leq T$. Especially, for non-path-dependent European options one can just solve the forward problem and approximate the non-discounted fair value at the spot by $\overline{P}_{LV}(T)^{\rm{T}}U_{LV}(0)$.

Now, consider a non-path-dependent option whose payoff is only dependent on the exchange rate $S_{T}$. Then 
$$ U_{LV,i}(0) = \boldsymbol{U}_{i,j}(0) $$
whenever $1\leq i \leq m_{1}, 1 \leq j \leq m_{2}$. 
It is readily verified that the semidiscretizations of the LV model and the SLV model define the same value at the spot, i.e.\ $U_{LV,i_{0}}(T) = \boldsymbol{U}_{i_{0},j_{0}}(T)$, if 
\begin{equation*}
\overline{\boldsymbol{P}}(T) e_v  = \overline{P}_{LV}(T).
\end{equation*}
This property is desirable for every possible maturity. Hence, one would like to have
\begin{equation}
\overline{\boldsymbol{P}}(\tau) e_v = \overline{P}_{LV}(\tau)
\label{eq:DiscreteProbabilitiesMatch}
\end{equation}
for all $\tau \geq 0$.
Notice that \eqref{eq:DiscreteProbabilitiesMatch} is equivalent to 
$$ \sum_{j=1}^{m_{2}} \boldsymbol{P}_{i,j}(\tau) w_{v,j} = P_{LV,i}(\tau) \qquad \mathrm{for} \ 1 \leq i \leq m_{1},\, \tau \geq 0,$$
which can be viewed as a semidiscrete analogue of \eqref{eq:SamePriceContinuous}.
Since $\overline{\boldsymbol{P}}(0) e_v = \overline{P}_{LV}(0)$, the condition \eqref{eq:DiscreteProbabilitiesMatch} is satisfied if 
\begin{equation}
\overline{\boldsymbol{P}}'(\tau) e_v = \overline{P}'_{LV}(\tau)
\label{eq:MatchModels}
\end{equation}
for all $ \tau > 0$.
From \eqref{eq:Pmatrix}, \eqref{eq:DiscretizationBCV} we directly obtain
\begin{equation*}
\overline{\boldsymbol{P}}'(\tau) e_v = \tfrac{1}{2} D^{\rm{T}}_{xx} L^{2}(\tau) \overline{\boldsymbol{P}}(\tau)\psi^{2}(\Lambda) e_v + (r_{d}-r_{f})D^{\rm{T}}_{x}\, \overline{\boldsymbol{P}}(\tau) e_v - \tfrac{1}{2} D^{\rm{T}}_{x} L^{2}(\tau) \overline{\boldsymbol{P}}(\tau)\psi^{2}(\Lambda) e_v.
\end{equation*}
If the (initially unspecified) diagonal matrix $L(\tau)$ is now defined through
\begin{equation}
L^{2}(\tau) \overline{\boldsymbol{P}}(\tau)\psi^{2}(\Lambda) e_v = L^{2}_{LV}(\tau)\overline{\boldsymbol{P}}(\tau) e_v,
\label{eq:LtoMatch}
\end{equation}
then 
\begin{equation*}
\overline{\boldsymbol{P}}'(\tau) e_v = \tfrac{1}{2} D^{\rm{T}}_{xx} L^{2}_{LV}(\tau) \overline{\boldsymbol{P}}(\tau) e_v + (r_{d}-r_{f})D^{\rm{T}}_{x}\, \overline{\boldsymbol{P}}(\tau) e_v - \tfrac{1}{2} D^{\rm{T}}_{x} L^{2}_{LV}(\tau)\overline{\boldsymbol{P}}(\tau) e_v.
\end{equation*}
Hence, it follows that \eqref{eq:MatchModels} holds whenever equation \eqref{eq:SemidiscretizationForwardEquationLV} has a unique solution.

Remark that in the definition \eqref{eq:LtoMatch} for the semidiscrete leverage function it is tacitly assumed that both vectors
$$ \overline{\boldsymbol{P}}(\tau)\psi^{2}(\Lambda) e_v \quad \mathrm{and} \quad \overline{\boldsymbol{P}}(\tau) e_v $$
only contain strictly positive values.
By performing a spatial discretization with finite differences it is possible that some of the values $\overline{\boldsymbol{P}}_{i,j}$ become negative. In our experiments, both vectors remained strictly positive, however, for natural values of $m_{1}, m_{2}$.

Notice that \textit{the derivation above is not restricted to the choice of finite difference formulas}.
If the second-order central formulas from Subsection \ref{subsec:BackwardDiscretization} are replaced by alternative finite difference formulas for which \eqref{eq:DiscretizationBCV} holds, and if these formulas are also applied for a similar semidiscretization in the LV model, then the SLV model is calibrated exactly to the LV model by employing \eqref{eq:LtoMatch}. 

We arrive at the following main result.

\begin{theorem}
Assume semidiscretization of the backward Kolmogorov equation \eqref{eq:BackwardKolmogorov} is performed by consistent finite difference formulas on a Cartesian grid and that semidiscretization of the forward Kolmogorov equation \eqref{eq:ForwardKolmogorov} is performed by the adjoint spatial discretization. Then
$$ U_{k_{0}}(T) = \overline{P}(\tau)^{\mathrm{T}}U(T-\tau) \qquad \mathit{for~all}\, \ 0 \le \tau \le T, $$
where $k_{0}$ corresponds to the index $(i_{0},j_{0})$ such that $(x_{i_{0}},v_{j_{0}}) = (X_{0},V_{0})$.
Hence, the two semidiscretizations define the same approximation for the fair value of non-path-dependent European options.

Next, assume semidiscretization of the backward and forward equations \eqref{eq:BackwardKolmogorovLV} and \eqref{eq:ForwardKolmogorovLV} under the LV model is performed in complete correspondence to that of \eqref{eq:BackwardKolmogorov} and \eqref{eq:ForwardKolmogorov}, respectively, under the SLV model by using the same grid and finite difference formulas in the $x$-direction.
If \eqref{eq:DiscretizationBCV} holds and equation \eqref{eq:SemidiscretizationForwardEquationLV} has a unique solution and if the leverage function $\sigma_{SLV}$ is defined on the grid in the $x$-direction by
\begin{equation}
\sigma^{2}_{SLV}(x_{i},\tau) =  \sigma^{2}_{LV}(x_{i},\tau) \frac{\sum_{j=1}^{m_{2}} \overline{\boldsymbol{P}}_{i,j}(\tau)}{\sum_{j=1}^{m_{2}}\psi^{2}(v_{j}) \overline{\boldsymbol{P}}_{i,j}(\tau)},
\label{eq:SigmaDiscretetoMatch}
\end{equation}
then
$$ \overline{\boldsymbol{P}}(\tau) e_v = \overline{P}_{LV}(\tau) \qquad \mathit{for~all}\, \ \tau \ge 0.  $$
In particular, if the payoff depends only on the exchange rate $S_{T}$, then the semidiscretizations of the LV model and the SLV model define the same approximation for the fair value of non-path-dependent European options:
$$ U_{LV,i_{0}}(T) = \overline{P}_{LV}(\tau)^{\mathrm{T}}U_{LV}(T-\tau) = \overline{P}(\tau)^{\mathrm{T}}U(T-\tau) = U_{k_{0}}(T) \qquad \mathit{for~all}\, \ 0 \le \tau \le T. $$
\label{th:LtoMatchTheorem}
\end{theorem}
The second part of Theorem \ref{th:LtoMatchTheorem} can be regarded as the semidiscrete analogue of \eqref{eq:SigmatoMatch}.
Indeed, if $\sigma_{SLV}$ is defined on the spatial grid in the $x$-direction by \eqref{eq:SigmaDiscretetoMatch}, then by the definition of $\overline{\boldsymbol{P}}$ it is directly seen that this is equivalent to applying \eqref{eq:SigmatoMatch}, where the conditional expectation is approximated by
\begin{equation}
\mathbb{E}[ \psi^{2}(V_{\tau}) \vert X_{\tau} = x_{i} ] \approx \frac{\sum_{j=1}^{m_{2}}\psi^{2}(v_{j}) \boldsymbol{P}_{i,j}(\tau) w_{v,j}}{\sum_{j=1}^{m_{2}} \boldsymbol{P}_{i,j}(\tau)w_{v,j}} = \frac{\sum_{j=1}^{m_{2}}\psi^{2}(v_{j}) \overline{\boldsymbol{P}}_{i,j}(\tau)}{\sum_{j=1}^{m_{2}} \overline{\boldsymbol{P}}_{i,j}(\tau)} 
\qquad \mathrm{for}~ 1 \leq i \leq m_{1}.
\label{eq:SemidiscreteExpectation}
\end{equation}

\setcounter{equation}{0}
\section{Time discretization}
\label{TimeDiscretization}

In this section we consider a suitable time stepping method for the numerical solution of semidiscrete systems of the type \eqref{eq:SemidiscretizationBackwardEquation} and \eqref{eq:SemidiscretizationForwardEquation} assuming that the matrix-valued function $L$ is known.
In general, semidiscretization by means of finite difference or finite volume methods of initial-boundary value problems for two-dimensional time-dependent convection-diffusion equations leads to large systems of stiff ODEs,
\begin{equation*}
W'(t)=F(t,W(t)) \quad (0 < t \leq T), \quad W(0)=W_{0},
\end{equation*}
with given vector-valued function $F: [0,T]\times \R^m \rightarrow \R^m$ and given vector $W_{0}\in \R^m$ where \mbox{$m\ge 1$} is the number of spatial grid points.
For the effective time discretization of these semidiscrete systems, \textit{operator splitting schemes} of the \textit{Alternating Direction Implicit} (ADI) type are widely used in practice. 
Several ADI schemes have been developed and analyzed in the literature for the situation where mixed spatial derivative terms are present in the convection-diffusion equation.
Mixed derivative terms are ubiquitous in the field of financial option valuation; they arise due to correlations between the underlying stochastic processes.
Notably, this holds for the general SLV model \eqref{eq:SLVmodelX} under consideration. 

In this paper we employ the \textit{Modified Craig--Sneyd (MCS) scheme}, which is a modern scheme of the ADI-type suitable to deal with mixed derivative terms.
Let the vector-valued function $F$ be decomposed as
\begin{equation*}
F(t,w) = F_{0}(t,w) + F_{1}(t,w) + F_{2}(t,w) \quad (0 \leq t \leq T,~w \in \R^m),
\end{equation*}
where $F_{0}$ represents the mixed spatial derivative term and $F_{1}$, respectively $F_{2}$, represents all spatial derivative terms in the first, respectively second, spatial direction.
Let $\theta > 0$ be a given parameter, $N \geq 1$ the number of time steps and $t_{n} = n \Delta t$ with $\Delta t = T/N$. 
Then the MCS scheme defines approximations $W_{n}$ to $W(t_{n})$ successively for $n=1,2,3,\ldots,N$ through
\begin{equation}
\label{eq:MCS}
\left\{\begin{array}{l}
Y_0 = W_{n-1}+\Delta t\, F(t_{n-1},W_{n-1}), \\\\
Y_l = Y_{l-1}+\theta\Delta t \left(F_l(t_n,Y_l)-F_l(t_{n-1},W_{n-1})\right),
\quad l=1,2, \\\\
\widehat{Y}_0 = Y_0+ \theta \Delta t \left(F_0(t_n,Y_2)-F_0(t_{n-1},W_{n-1})\right),\\\\
\widetilde{Y}_0 = \widehat{Y}_0+ (\tfrac{1}{2}-\theta )\Delta t \left(F(t_n,Y_2)-F(t_{n-1},W_{n-1})\right),\\\\
\widetilde{Y}_l = \widetilde{Y}_{l-1}+\theta\Delta t \,(F_l(t_n,\widetilde{Y}_l)-F_l(t_{n-1},W_{n-1})),
\quad l=1,2, \\\\
W_n = \widetilde{Y}_2.
\end{array}\right.
\end{equation}
The MCS scheme \eqref{eq:MCS} has been introduced by in 't Hout \& Welfert \cite{IHW09}.
It starts with an explicit Euler predictor stage, which is followed by two implicit but unidirectional corrector stages.
Then a second explicit stage is performed, followed again by two implicit unidirectional corrector stages.
In the special case where $\theta = \frac{1}{2}$ one obtains the {\it Craig--Sneyd (CS) scheme}, proposed in \cite{CS88}.
Observe the favourable feature that the $F_0$ part, representing the mixed derivative term, is always treated in an explicit manner.

In the past years various positive \textit{stability} results have been derived for the MCS scheme relevant to multidimensional convection-diffusion equations with mixed derivative terms, see e.g.~\cite{IHM11,IHM13,IHW09,M14}. Subsequently, in 't Hout \& Wyns~\cite{IHW15} proved that, under natural stability and smoothness assumptions, the MCS scheme is \textit{second-order convergent} with respect to the time step whenever it is applied to semidiscrete two-dimensional convection-diffusion equations with mixed derivative term. The temporal convergence bound from \cite{IHW15} has the key property that it holds uniformly in the spatial mesh width.

In the present application, both vectors $U(0)$ and $\overline{P}(0)$ are stemming from initial functions that are nonsmooth. 
It is well-known in the literature that the convergence behaviour of time discretization methods can be seriously impaired for nonsmooth initial data, compare~\cite{PVF03,W16}.
To alleviate this, Rannacher time stepping will be applied, that is, the first few time steps of the MCS scheme are replaced by twice as many half time steps using the implicit Euler scheme \cite{R84}.
We opt to replace here the first two MCS time steps by four half time steps of the implicit Euler scheme. 
Based on the recent results in \cite{W16}, the MCS scheme is then expected to maintain second-order convergence in time, uniformly in the number of spatial grid points.

As seen in Section~\ref{ForwardBackward}, the fair values of non-path-dependent European options can be determined by solving either the backward equation \eqref{eq:BackwardKolmogorov} or the forward equation \eqref{eq:ForwardKolmogorov}. 
In Section~\ref{BackwardForward} spatial discretization of these two PDEs led to the semidiscrete systems \eqref{eq:SemidiscretizationBackwardEquation} and \eqref{eq:SemidiscretizationForwardEquation}, respectively.
Application of the MCS scheme to \eqref{eq:SemidiscretizationBackwardEquation} yields approximations $U_{n}$ of $U(t_{n})$ and the non-discounted fair option value at the spot is then approximated by $U_{k_{0},N} = \overline{P}^{\rm{T}}_{0} U_{N}$.
Alternatively, take $\Delta \tau = \Delta t = T/N$ and let temporal grid points $\tau_{n} = n \Delta \tau = T - t_{N-n}$.
Application of the MCS scheme to \eqref{eq:SemidiscretizationForwardEquation} yields approximations $\overline{P}_{n}$ of $\overline{P}(\tau_{n})$ and the non-discounted fair option value at the spot is then approximated by $\overline{P}_{N}^{\rm{T}}U_0$.

It is possible, see Itkin~\cite{I15}, to construct a new ADI discretization for \eqref{eq:SemidiscretizationForwardEquation} such that there is an exact match between the fully discretized backward equation and the fully discretized forward equation, that is,
$$ \overline{P}^{\rm{T}}_{0} U_{N} = \overline{P}^{\rm{T}}_{N} U_{0} = \overline{P}^{\rm{T}}_{N-n} U_{n} $$
for all $0 \leq n \leq N$.
In this paper we prefer to employ the MCS scheme for the numerical solution of both ODE systems \eqref{eq:SemidiscretizationBackwardEquation} and \eqref{eq:SemidiscretizationForwardEquation}.
A main reason is that ample positive results are already available in the literature on the stability and convergence of the MCS scheme.
The alternative scheme from \cite{I15} ends with an explicit stage and numerical experiments indicate that computing \eqref{eq:SemidiscreteExpectation} with the corresponding $P_{n}$ can lead to undesirable (erratic) behaviour of the numerical conditional expectation.
In addition, practical experience shows that the temporal discretization error of the MCS scheme is often much smaller than the spatial discretization error.

\setcounter{equation}{0}
\section{Calibration of the SLV model to the LV model}\label{Calibration}

In Section \ref{LocalVolMatching} we derived the expression \eqref{eq:SigmaDiscretetoMatch} for the discrete leverage function $\sigma_{SLV}$ that exactly calibrates the semidiscrete SLV model to the semidiscrete LV model.
This expression involves the matrix function $\overline{\boldsymbol{P}}$.
Combining \eqref{eq:SigmaDiscretetoMatch} with the semidiscrete forward equation \eqref{eq:SemidiscretizationForwardEquation} for $\overline{P} = \mathrm{vec}[\overline{\boldsymbol{P}}]$, one arrives at a large, non-linear system of ODEs.
In the present section numerical time-stepping is applied together with an inner iteration so as to numerically solve this system of ODEs and acquire the discrete leverage function that satisfies \eqref{eq:SigmaDiscretetoMatch}.

Suppose an approximation $\overline{P}_{n}$ to $\overline{P}(\tau_{n})$ at time level $\tau_{n}$ is known.
Let $\overline{\boldsymbol{P}}_{n}$ denote the $m_{1} \times m_{2}$ matrix such that $$ \overline{P}_{n} = \mathrm{vec}[\overline{\boldsymbol{P}}_{n}]. $$
Then the discrete leverage function is determined by 
\begin{equation}
\sigma^{2}_{SLV}(x_{i},\tau_{n}) \mathbb{E}_{n,i} = \sigma^{2}_{LV}(x_{i},\tau_{n}),
\label{eq:FullyDiscreteLeverage}
\end{equation}
where the quantity
\begin{equation}
\mathbb{E}_{n,i} = \frac{\sum_{j=1}^{m_{2}}\psi^{2}(v_{j}) \overline{\boldsymbol{P}}_{n,i,j}}{\sum_{j=1}^{m_{2}} \overline{\boldsymbol{P}}_{n,i,j}}
\label{eq:FullyDiscreteExpectation}
\end{equation}
forms an approximation to the conditional expectation $\mathbb{E}[ \psi^{2}(V_{\tau_{n}}) \vert X_{\tau_{n}} = x_{i} ]$, see \eqref{eq:SemidiscreteExpectation}.
In order to arrive at the actual calibration procedure, we need to consider three practical issues concerning formula \eqref{eq:FullyDiscreteExpectation}.

$\bullet$~ Due to the spatial and temporal discretizations, it may happen that either the numerator or denominator of \eqref{eq:FullyDiscreteExpectation} becomes negative. 
In this case we assume that the conditional expectation is locally constant in time and set $\mathbb{E}_{n,i} = \mathbb{E}_{n-1,i}$. 
In our experiments, however, both parts of the quotient remained strictly positive for common values of $m_{1}, m_{2}, \Delta \tau$.

$\bullet$~ Since $\overline{\boldsymbol{P}}_{n,i,j}$ is an approximation of the probability of the event that
$$ (X_{\tau_{n}}, V_{\tau_{n}}) \in [x_{i} - \tfrac{\Delta x_{i}}{2}, x_{i} + \tfrac{\Delta x_{i+1}}{2}] \times [v_{j} - \tfrac{\Delta v_{j}}{2}, v_{j} + \tfrac{\Delta v_{j+1}}{2}], $$
it can happen that the numerator and denominator of \eqref{eq:FullyDiscreteExpectation} become very small, which can lead to unrealistic values of the leverage function.
To resolve this, a regularized approximation of the conditional expectation is used (cf.~\cite{EKO12}),
\begin{equation}
\mathbb{E}_{n,i} = \frac{\sum_{j=1}^{m_{2}}\psi^{2}(v_{j}) \overline{\boldsymbol{P}}_{n,i,j} + \psi^{2}(\eta) \epsilon}{\sum_{j=1}^{m_{2}} \overline{\boldsymbol{P}}_{n,i,j} + \epsilon}
\label{eq:FullyDiscreteExpectationStable}
\end{equation}
for given small value $\epsilon$. In this paper $\epsilon = 10^{-8}$ is taken. By using the regularized version \eqref{eq:FullyDiscreteExpectationStable}, the approximated conditional expectation is shifted towards $\psi^{2}(\eta)$ where $\eta$ is the mean-reversion level of the process $V_{\tau}$.

$\bullet$~ At the spot $\tau=0$, the matrix $\overline{\boldsymbol{P}}(0)$ has $(i_{0},j_{0})$-th entry equal to one and all its other entries are equal to zero.
Consequently, if $n=0$, then the expression \eqref{eq:FullyDiscreteExpectation} is only defined if $i=i_{0}$.
To render the calibration procedure feasible, we extend this definition to all indices $i$ and thus put
$$ \mathbb{E}_{0,i} = \psi^{2}(V_{0}). $$
Notice that this agrees with \eqref{eq:FullyDiscreteExpectationStable} for $n=0$ whenever $\eta = V_{0}$, which often holds in practice.\\

Let $Q\ge 1$ be a given integer.
For calibrating the SLV model to the LV model, we employ the following numerical procedure.
It consists of numerical time stepping combined with an inner iteration, cf.~\cite{TF10}.
\newline
\newline
\texttt{for $n$ is $1$ to $N$ do
\begin{itemize}
\item[] let $\overline{P}_{n} = \overline{P}_{n-1}$ be an initial approximation to $\overline{P}(\tau_{n})$; \newline \newline 
for $q$ is $1$ to $Q$ do
\begin{itemize}
\item[(a)] approximate the conditional expectations $\mathbb{E}[ \psi^{2}(V_{\tau_{n}}) \vert X_{\tau_{n}} = x_{i} ]$ by \eqref{eq:FullyDiscreteExpectationStable};
\item[(b)] approximate $\sigma_{SLV}(\cdot,\tau_{n})$ on the grid in the $x$-direction by formula \eqref{eq:FullyDiscreteLeverage};
\item[(c)] update $\overline{P}_{n}$ by performing a numerical time step for \eqref{eq:SemidiscretizationForwardEquation} from $\tau_{n-1}$ to $\tau_{n}$;
\end{itemize}
end
\end{itemize}
end  
\newline \newline
}
Whenever a time step from $\tau_{n-1}$ to $\tau_{n}$ with the MCS scheme is replaced by two half-time steps of the implicit Euler scheme, the inner iteration above is first performed for the substep from $\tau_{n-1}$ to $\tau_{n-1/2} = \tau_{n-1}+\Delta \tau /2$, yielding an approximation of $\overline{P}(\tau_{n-1/2})$ and $\sigma_{SLV}(\cdot,\tau_{n-1/2})$. Next, the inner iteration is performed for the substep from  $\tau_{n-1/2}$ to $\tau_{n}$, yielding an approximation of $\overline{P}(\tau_{n})$ and $\sigma_{SLV}(\cdot,\tau_{n})$.

Upon completion of the time stepping and iteration procedure above, the original approximation for $\sigma_{SLV}(\cdot,0)$ is replaced on the grid in the $x$-direction by
$$ \sigma^{2}_{LV}(x_{i},0) = \sigma^{2}_{SLV}(x_{i},0) \frac{\sum_{j=1}^{m_{2}}\psi^{2}(v_{j}) \overline{\boldsymbol{P}}_{i,j,1} + \psi^{2}(\eta) \epsilon}{\sum_{j=1}^{m_{2}} \overline{\boldsymbol{P}}_{i,j,1} + \epsilon} $$
for $1 \leq i \leq m_{1}$. This appears more realistic as the original approximation was actually only valid for the index $i=i_{0}$.

\setcounter{equation}{0}
\section{Numerical experiments}\label{Numerical experiments}

In this section, numerical experiments are presented to illustrate the effectiveness of the calibration procedure.
Here, we opt to consider the popular and challenging Heston-based SLV model, i.e.\ SLV model \eqref{eq:SLVmodelX} with $\psi(v) = \sqrt{v}$ and $\alpha = 1/2$, to describe the evolution of the EUR/USD exchange rate. 

As stated in the introduction, for the calibration of the SLV model it is customary to start from the parameters $\kappa,\eta,\xi_{SV},\rho$ of the underlying (Heston) SV model such that this model reflects the market dynamics of the exchange rate.  
Next a mixing parameter $0 \le \mu \le 1$ is used to define $\xi = \mu \xi_{SV}$ and thus to tune between the underlying LV and SV models. 
In this article we consider the following four sets of parameters:
\begin{center}
\begin{tabular}{|c| r r r r r r r|}
\hline \newline
& $\kappa$ & $\eta$ & $\xi_{SV}$ & $\rho$ & $\mu$ & $\xi$ & $T$ \\
\hline \newline
Case 1 & $3.02$ & $ 0.015 $ & $ 0.41 $ & $ -0.13 $ & $ 0.75 $ & $ 0.31 $ &  $ 6M $  \\
\hline
\newline
Case 2 & $ 1.00 $ & $ 0.09 $ & $ 1.00 $ & $ -0.3 $ & $ 1 $ & $ 1.00 $ &  $ 6M $ \\
\hline
\newline
Case 3 & $ 0.75 $ & $ 0.015 $ & $ 0.20 $ & $ -0.14 $ & $ 0.75 $ & $ 0.15 $ &  $ 2Y $ \\
\hline
\newline
Case 4 & $ 1.00 $ & $ 0.09 $ & $ 1.00 $ & $ -0.3 $ & $ 1 $ & $ 1.00 $ &  $ 2Y $ \\
\hline
\end{tabular}
\end{center}
The first and third parameter set are taken from \cite{C11}. They correspond to the EUR/USD exchange rate for the pertinent maturities (market data as of 16 September 2008). 
The second and fourth parameter set are essentially the same and taken from \cite{A08}. The latter two sets are challenging for the calibration procedure as the \textit{Feller condition} is strongly violated, i.e.\ $2\kappa \eta \ll \sigma^{2}$, so that the probability mass is stacked up near $v=0$.

\begin{figure}
\begin{center}
\includegraphics[scale=0.5]{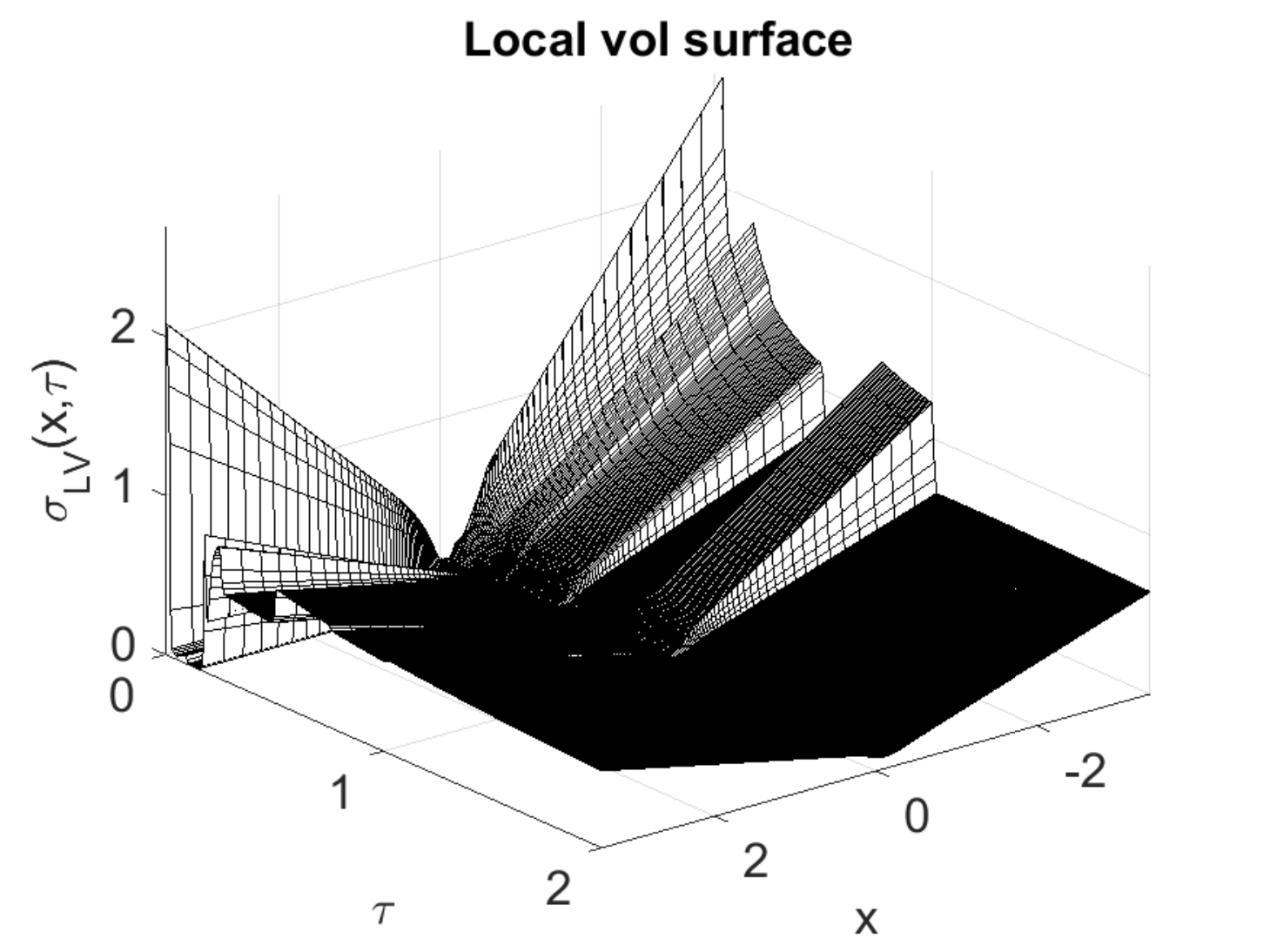} 
\includegraphics[scale=0.5]{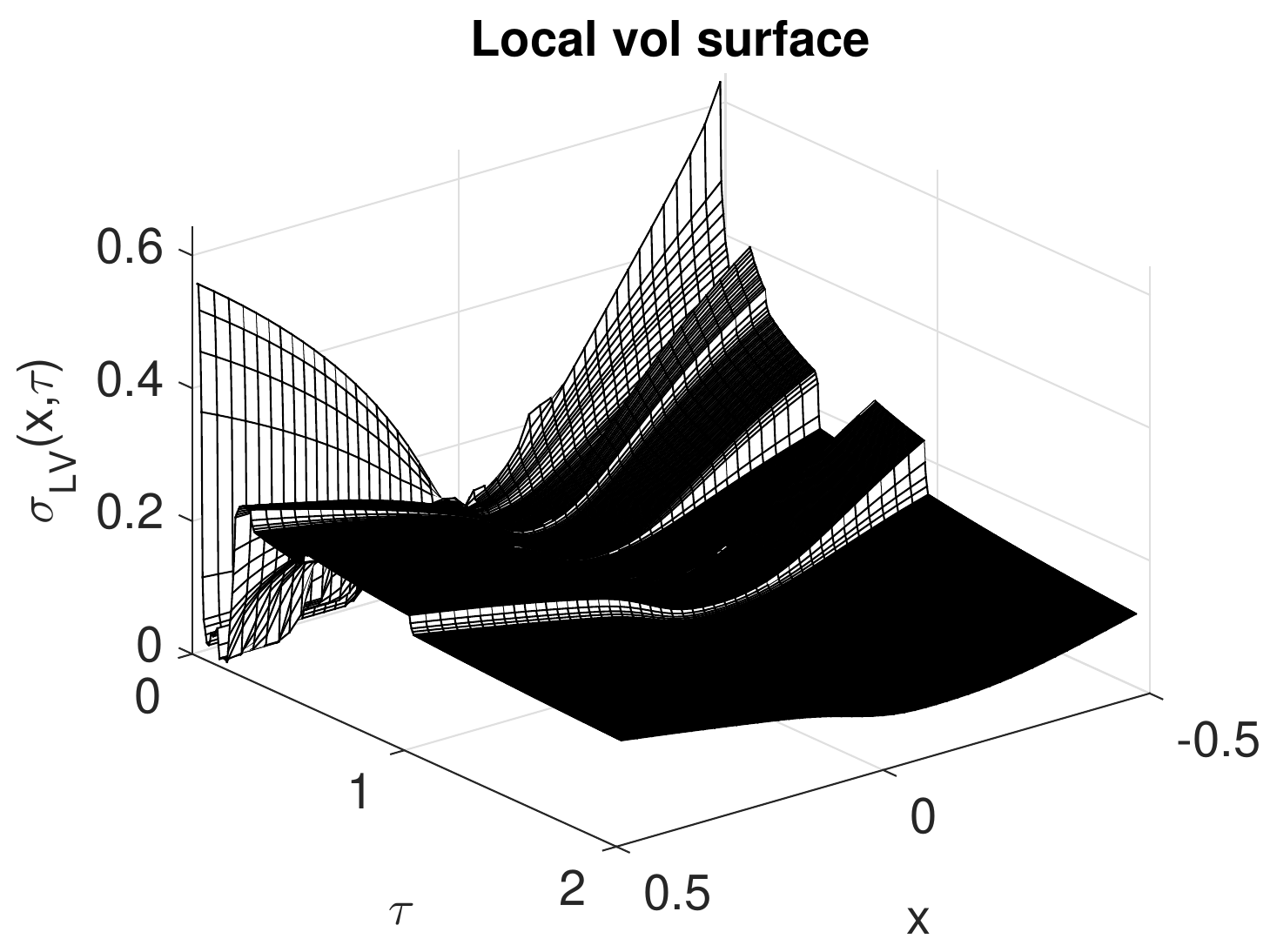} 
\caption{Local volatility function originating from actual EUR/USD vanilla option data (market data as of 13 November 2015).
The spot rate $S_{0} = 1.0764$.}
\label{fig:LVSurface}
\end{center}
\end{figure}

The leverage function $\sigma_{SLV}$ is determined in such a way that the SLV model is calibrated exactly to the underlying LV model. The LV model is completely determined by the LV function $\sigma_{LV}$ and the risk-free interest rates $r_{d}, r_{f}$. For the experiments we consider
$$ r_{d} = 0.03, \quad r_{f} = 0.01, $$
and LV function displayed in Figure \ref{fig:LVSurface}.
This LV function originates from actual EUR/USD vanilla option data (market data as of 13 November 2015). 
The corresponding spot rate is
$$ S_{0} = 1.0764. $$
For the spot value $V_{0}$ of the process $V_{\tau}$ in the SLV model we assume that it is equal to the long-term mean $\eta$ of this process.

\begin{figure}
\begin{center}
\includegraphics[scale=0.5]{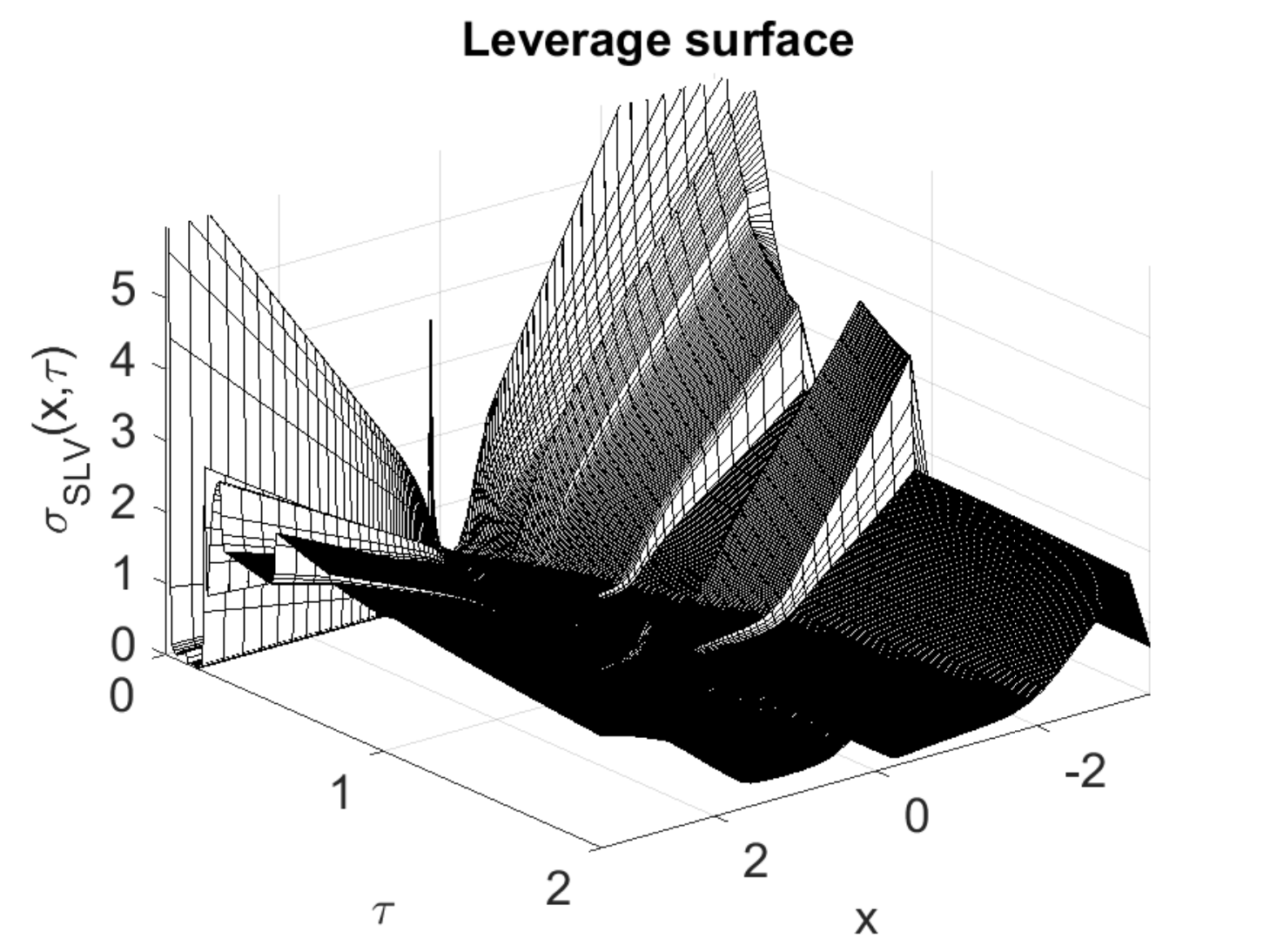} 
\includegraphics[scale=0.5]{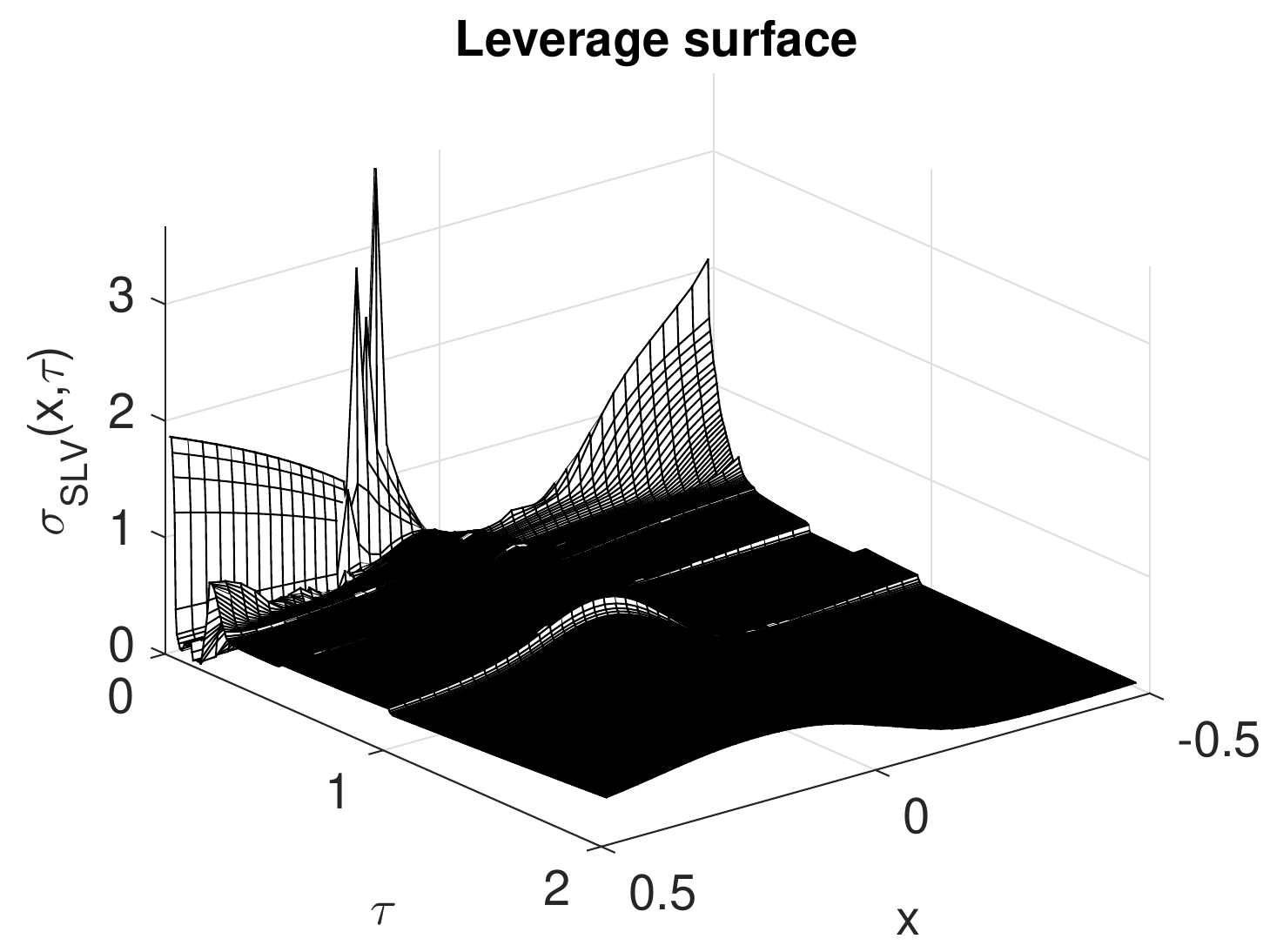} 
\caption{Leverage function stemming from the calibration procedure with local volatility function from Figure \ref{fig:LVSurface}, SV parameters from Case 4 and values $m_{1} = 100$, $m_{2} = 50$, $\Delta \tau = 1/200$, $\theta = 1/3$, $Q =2.$}
\label{fig:SLVSurface}
\end{center}
\end{figure}

The aim is to construct the leverage function in such a way that the discretized LV model and the discretized SLV model yield identical approximate values for any given vanilla option whenever similar discretizations are employed. 
In our experiments, the semidiscretization of the backward Kolmogorov equation \eqref{eq:BackwardKolmogorov} is performed as described in Subsection \ref{subsec:BackwardDiscretization} and the semidiscretization of the forward equation \eqref{eq:ForwardKolmogorov} is defined by the pertinent adjoint spatial discretization. The backward and forward PDEs \eqref{eq:BackwardKolmogorovLV} and \eqref{eq:ForwardKolmogorovLV} are semidiscretized analogously and by using the same finite difference schemes as described in Section \ref{LocalVolMatching}.
For the first numerical experiment we consider
$$ m_{1} = 100, \qquad m_{2} = 50. $$
The main Theorem \ref{th:LtoMatchTheorem} yields that if the leverage function is defined on the grid in the $x$-direction by \eqref{eq:SigmaDiscretetoMatch}, then the approximations obtained from the four semidiscrete systems \eqref{eq:SemidiscretizationBackwardEquation}, \eqref{eq:SemidiscretizationForwardEquation}, \eqref{eq:SemidiscretizationBackwardEquationLV},~\eqref{eq:SemidiscretizationForwardEquationLV} of the fair value of any given non-path-dependent European option are identical.
The exact solution \eqref{eq:SigmaDiscretetoMatch} is approximated by applying the calibration procedure described in Section \ref{Calibration}. We choose to perform the temporal discretization in this procedure with values
$$ \Delta \tau = 1/200, \qquad \theta = 1/3, \qquad Q = 2.$$
In Figure \ref{fig:SLVSurface} the obtained discrete leverage function is shown for Case 4.
If the SV model with parameters from Case 4 would fit the market prices for European call and put options exactly, then the leverage function would be identically equal to one and the SLV model reduces to the SV model. Clearly, Figure \ref{fig:SLVSurface} indicates that the pure SV model with parameters from Case 4 doesn't match the market data very well. 
This outcome was to be expected, as the SV parameters from \cite{A08} do not correspond to a EUR/USD exchange rate.

With the obtained discrete leverage function, the performance of the calibration procedure is first tested by comparing the fully discrete approximations of the fair value of European call options which are acquired by numerically solving the systems of ODEs \eqref{eq:SemidiscretizationBackwardEquation}, \eqref{eq:SemidiscretizationForwardEquation}, \eqref{eq:SemidiscretizationBackwardEquationLV}, \eqref{eq:SemidiscretizationForwardEquationLV}.
To this purpose we consider a range of strikes, given by
$$ K = 0.7 S_{0}, \ 0.8 S_{0}, \ 0.9 S_{0}, \ S_{0}, \ 1.1 S_{0}, \ 1.2 S_{0}, \ 1.3 S_{0}. $$
Temporal discretization of \eqref{eq:SemidiscretizationBackwardEquationLV} and \eqref{eq:SemidiscretizationForwardEquationLV} is performed by the classical \textit{Crank--Nicolson scheme}. 
The systems \eqref{eq:SemidiscretizationBackwardEquation} and \eqref{eq:SemidiscretizationForwardEquation}, which are stemming from a two-dimensional PDE, are discretized in time by the MCS scheme \eqref{eq:MCS} with parameter $\theta$ given above.
For both methods we consider $\Delta \tau = 1/200$ and Rannacher time stepping is applied to handle the nonsmoothness of the initial functions.
In Table \ref{Table:Case1m100N200} the obtained fully discrete approximations $(FV_{m})$, $m \in \{LVB, LVF, SLVB, SLVF\}$, of the fair value (FV) are presented in Case 1. 
Here $m=LVB$, respectively $m=LVF,SLVB,SLVF$, corresponds with the approximated fair value obtained via \eqref{eq:SemidiscretizationBackwardEquationLV}, respectively \eqref{eq:SemidiscretizationForwardEquationLV}, \eqref{eq:SemidiscretizationBackwardEquation}, \eqref{eq:SemidiscretizationForwardEquation}. 
In Table \ref{Table:Case1m100N200} one observes that the approximated option values are almost identical. To express in more detail the quality of the approximations, we  present relative errors
$$\epsilon_{r,m} = (FV_{m}-FV_{LVB})/FV_{LVB}.$$
Here the option values given by solving \eqref{eq:SemidiscretizationBackwardEquationLV}, indicated by $FV_{LVB}$, are considered as the reference values. This is motivated by the fact that in practice one starts from the underlying LV model and within the LV model it is common to solve the backward equation \eqref{eq:BackwardKolmogorovLV}.
Table \ref{Table:Case1m100N200} reveals the favourable result that all relative errors are smaller than $0.1\%$. 
Numerical experiments for the other SV parameter sets, i.e.\ for Cases $2,3,4$, yield the same observation.
It can be concluded that the different approximations are almost identical in each of the four cases and the calibration procedure from Section \ref{Calibration} performs well. 

\begin{table}
\begin{center}
\begin{tabular}{|c|c||c|c||c|c||c|c|}
\hline
$K/S_{0}$ & $FV_{LVB}$ & $FV_{LVF}$ & $\epsilon_{r,LVF}$ & $FV_{SLVB}$ & $\epsilon_{r,SLVB}$ & $FV_{SLVF}$ & $\epsilon_{r,SLVF}$ \\
\hline
$0.7$ & $0.3288$ & $0.3288$ & $0.0000 \%$ & $0.3288$ & $0.0000 \%$ & $0.3288$ & $0.0000 \%$ \\
$0.8$ & $0.2228$ & $0.2228$ & $0.0000 \%$ & $0.2228$ & $0.0000 \%$ & $0.2228$ & $0.0001 \%$ \\
$0.9$ & $0.1185$ & $0.1185$ & $0.0008 \%$ & $0.1185$ & $0.0004 \%$ & $0.1185$ & $0.0008 \%$ \\
$1$   & $0.0381$ & $0.0381$ & $0.0083 \%$ & $0.0381$ & $0.0004 \%$ & $0.0381$ & $0.0079 \%$ \\
$1.1$ & $0.0091$ & $0.0091$ & $0.0235 \%$ & $0.0091$ & $0.0017 \%$ & $0.0091$ & $0.0211 \%$ \\
$1.2$ & $0.0019$ & $0.0019$ & $0.0437 \%$ & $0.0019$ & $0.0073 \%$ & $0.0019$ & $0.0391 \%$ \\
$1.3$ & $0.0004$ & $0.0004$ & $0.0661 \%$ & $0.0004$ & $0.0223 \%$ & $0.0004$ & $0.0640 \%$ \\
\hline
\end{tabular}
\end{center}
\caption{Comparison of the approximated option values $FV_{LVB}, FV_{LVF}, FV_{SLVB}, FV_{SLVF}$ in Case 1 and for values $m_{1} = 100, m_{2} = 50, \Delta \tau = 1/200, \theta = 1/3, Q =2.$}
\label{Table:Case1m100N200}
\end{table}

When the strike increases relative to $S_0$, the fair value of European call options tends to zero and it is difficult to adequately compare approximations. In financial practice, European call and put options are often quoted in terms of {\it implied volatility}. Let $\sigma_{imp,m}$ denote the implied volatility (in~$\%$) corresponding to $FV_{m}$. In the following we test the performance of the calibration procedure by calculating the absolute implied volatility errors
$$ \epsilon_{m} = \vert \sigma_{imp,m} - \sigma_{imp,LVB} \vert . $$
In Table \ref{Table:ImpliedVolm100N200} these errors are presented for the four different SV parameter sets, taking the same values of $m_{1}, m_{2}, \Delta \tau, \theta, Q$ as above.
Since semidiscretization of \eqref{eq:ForwardKolmogorovLV} is performed by the adjoint spatial discretization, the only source of error in $\epsilon_{LVF}$ is the temporal discretization error. 
Table~\ref{Table:ImpliedVolm100N200} shows that these errors are small.
The somewhat larger values $\epsilon_{LVF}$ for $T=6M$ compared to $T=2Y$ can be explained from the fact that the implied volatility is more sensitive to changes in the fair value when the maturity is low. 
Table \ref{Table:ImpliedVolm100N200} subsequently reveals that in all numerical experiments the (small) errors $\epsilon_{SLVB}$ and $\epsilon_{SLVF}$ are of the same order of magnitude as $\epsilon_{LVF}$. This indicates that the size of the error due to the calibration (with the iteration to handle the non-linearity) is not larger than the size of the temporal discretization error. As the calibration procedure includes numerical time stepping, this is the best result one can aim for.

In order to verify the assertion above, we repeat the numerical experiments with a smaller step size. 
In Table \ref{Table:ImpliedVolm100N400} the absolute implied volatility errors are given for the same parameters as above but where the calibration and pricing is performed with 
$$ \Delta \tau = 1/400.$$
Comparing $\epsilon_{LVF}$ in Tables~\ref{Table:ImpliedVolm100N200} and~\ref{Table:ImpliedVolm100N400}, it is clearly seen that the temporal error decreases if $\Delta \tau$ decreases. Moreover, the absolute implied volatility errors $\epsilon_{SLVB}, \epsilon_{SLVF}$ are again of the same size as $\epsilon_{LVF}$. It can be concluded that for realistic values of $\Delta \tau$ the error introduced by the calibration procedure is of the same order of magnitude as the temporal discretization error.

Observe that by decreasing the step size, the reference values $\sigma_{imp,LVB}$ have slightly changed.
This is a consequence of the fact that more points from the LV surface are used and the fully discrete solution $FV_{LVB}$ converges to the semidiscrete solution $U_{LV,i_{0}}(T)$.
We note that the approximated option values $FV_{LVB}$ acquired with $\Delta \tau = 1/400$ are identical to those in Table~\ref{Table:Case1m100N200} up to the number of digits presented in that table. This confirms again that the temporal discretization error is small.

As stated in Theorem \ref{th:LtoMatchTheorem}, the condition \eqref{eq:SigmaDiscretetoMatch} facilitates an exact match between the semidiscrete LV model and the semidiscrete SLV model whenever similar discretizations are used. 
This match is valid for any number of spatial grid points $m_{1}, m_{2}$. In order to test this property of the calibration procedure, we repeat the numerical experiments with the number of spatial grid points replaced by
$$ m_{1} = 200, \quad m_{2} = 100, $$
and with step size $\Delta \tau = 1/200$.
The obtained implied volatilities $\sigma_{LVB}$ are presented in Table \ref{Table:ImpliedVolm200N200} as well as the absolute implied volatility errors $\epsilon_{LVF}, \epsilon_{SLVB}, \epsilon_{SLVF}$.

\begin{table}
\begin{center}
\begin{tabular}{|c|c|c||c|c||c|c|}
\hline
\multicolumn{3}{|c||}{$T = 6M$} & \multicolumn{2}{c||}{Case 1} & \multicolumn{2}{c|}{Case 2} \\
\hline
$K/S_{0}$ & $\sigma_{imp,LVB}$ & $\epsilon_{LVF}$ & $\epsilon_{SLVB}$ & $\epsilon_{SLVF}$ & $\epsilon_{SLVB}$ & $\epsilon_{SLVF}$ \\
\hline
$0.7$ & $14.8646$ & $0.0024$ & $0.0024$ & $0.0028$ & $0.0020$ & $0.0025$ \\
$0.8$ & $12.0458$ & $0.0017$ & $0.0015$ & $0.0018$ & $0.0014$ & $0.0015$ \\
$0.9$ & $10.3061$ & $0.0012$ & $0.0006$ & $0.0012$ & $0.0005$ & $0.0009$ \\
$1$   & $10.7970$ & $0.0011$ & $0.0000$ & $0.0010$ & $0.0006$ & $0.0007$ \\
$1.1$ & $12.6393$ & $0.0011$ & $0.0001$ & $0.0010$ & $0.0015$ & $0.0003$ \\
$1.2$ & $13.9656$ & $0.0011$ & $0.0002$ & $0.0010$ & $0.0002$ & $0.0005$ \\
$1.3$ & $14.9511$ & $0.0012$ & $0.0004$ & $0.0012$ & $0.0008$ & $0.0012$ \\
\hline
\end{tabular}
\end{center}
\begin{center}
\begin{tabular}{|c|c|c||c|c||c|c|}
\hline
\multicolumn{3}{|c||}{$T = 2Y$} & \multicolumn{2}{c||}{Case 3} & \multicolumn{2}{c|}{Case 4} \\
\hline
$K/S_{0}$ & $\sigma_{imp,LVB}$ & $\epsilon_{LVF}$ & $\epsilon_{SLVB}$ & $\epsilon_{SLVF}$ & $\epsilon_{SLVB}$ & $\epsilon_{SLVF}$ \\
\hline
$0.7$ & $10.2284$ & $0.0008$ & $0.0005$ & $0.0008$ & $0.0006$ & $0.0007$ \\
$0.8$ & $ 9.1864$ & $0.0005$ & $0.0003$ & $0.0005$ & $0.0003$ & $0.0004$ \\
$0.9$ & $ 8.9874$ & $0.0004$ & $0.0001$ & $0.0004$ & $0.0001$ & $0.0003$ \\
$1$   & $ 9.6063$ & $0.0004$ & $0.0000$ & $0.0004$ & $0.0001$ & $0.0003$ \\
$1.1$ & $10.6956$ & $0.0004$ & $0.0000$ & $0.0004$ & $0.0000$ & $0.0003$ \\
$1.2$ & $11.6810$ & $0.0004$ & $0.0001$ & $0.0004$ & $0.0000$ & $0.0002$ \\
$1.3$ & $12.4844$ & $0.0004$ & $0.0001$ & $0.0004$ & $0.0000$ & $0.0001$ \\
\hline
\end{tabular}
\end{center}
\caption{Comparison of the approximated implied volatilities $\sigma_{LVB}, \sigma_{LVF}, \sigma_{SLVB}, \sigma_{SLVF}$ for values $m_{1} = 100, m_{2} = 50, \Delta \tau = 1/200, \theta = 1/3, Q =2.$}
\label{Table:ImpliedVolm100N200}
\end{table}
\begin{table}
\begin{center}
\begin{tabular}{|c|c|c||c|c||c|c|}
\hline
\multicolumn{3}{|c||}{$T = 6M$} & \multicolumn{2}{c||}{Case 1} & \multicolumn{2}{c|}{Case 2} \\
\hline
$K/S_{0}$ & $\sigma_{imp,LVB}$ & $\epsilon_{LVF}$ & $\epsilon_{SLVB}$ & $\epsilon_{SLVF}$ & $\epsilon_{SLVB}$ & $\epsilon_{SLVF}$ \\
\hline
$0.7$ & $14.8605 $ & $0.0001 $ & $0.0007 $ & $0.0006 $ & $0.0006 $ & $0.0006 $ \\
$0.8$ & $12.0438 $ & $0.0000 $ & $0.0002 $ & $0.0001 $ & $0.0001 $ & $0.0000 $ \\
$0.9$ & $10.3053 $ & $0.0000 $ & $0.0000 $ & $0.0000 $ & $0.0001 $ & $0.0001 $ \\
$1$   & $10.7964 $ & $0.0000 $ & $0.0000 $ & $0.0000 $ & $0.0001 $ & $0.0001 $ \\
$1.1$ & $12.6386 $ & $0.0000 $ & $0.0000 $ & $0.0000 $ & $0.0006 $ & $0.0003 $ \\
$1.2$ & $13.9647 $ & $0.0000 $ & $0.0001 $ & $0.0000 $ & $0.0003 $ & $0.0001 $ \\
$1.3$ & $14.9498 $ & $0.0000 $ & $0.0001 $ & $0.0000 $ & $0.0001 $ & $0.0000 $ \\
\hline
\end{tabular}
\end{center}
\begin{center}
\begin{tabular}{|c|c|c||c|c||c|c|}
\hline
\multicolumn{3}{|c||}{$T = 2Y$} & \multicolumn{2}{c||}{Case 3} & \multicolumn{2}{c|}{Case 4} \\
\hline
$K/S_{0}$ & $\sigma_{imp,LVB}$ & $\epsilon_{LVF}$ & $\epsilon_{SLVB}$ & $\epsilon_{SLVF}$ & $\epsilon_{SLVB}$ & $\epsilon_{SLVF}$ \\
\hline
$0.7$ & $10.2278 $ & $0.0001 $ & $0.0001 $ & $0.0001 $ & $0.0001 $ & $0.0001 $ \\
$0.8$ & $ 9.1859 $ & $0.0000 $ & $0.0000 $ & $0.0000 $ & $0.0000 $ & $0.0000 $ \\
$0.9$ & $ 8.9870 $ & $0.0000 $ & $0.0000 $ & $0.0000 $ & $0.0000 $ & $0.0000 $ \\
$1$   & $ 9.6060 $ & $0.0000 $ & $0.0000 $ & $0.0000 $ & $0.0000 $ & $0.0000 $ \\
$1.1$ & $10.6953 $ & $0.0000 $ & $0.0000 $ & $0.0000 $ & $0.0000 $ & $0.0000 $ \\
$1.2$ & $11.6806 $ & $0.0000 $ & $0.0000 $ & $0.0000 $ & $0.0000 $ & $0.0000 $ \\
$1.3$ & $12.4840 $ & $0.0000 $ & $0.0000 $ & $0.0000 $ & $0.0001 $ & $0.0001 $ \\
\hline
\end{tabular}
\end{center}
\caption{Comparison of the approximated implied volatilities $\sigma_{LVB}, \sigma_{LVF}, \sigma_{SLVB}, \sigma_{SLVF}$ for values $m_{1} = 100, m_{2} = 50, \Delta \tau = 1/400, \theta = 1/3, Q =2.$}
\label{Table:ImpliedVolm100N400}
\end{table}

\begin{table}
\begin{center}
\begin{tabular}{|c|c|c||c|c||c|c|}
\hline
\multicolumn{3}{|c||}{$T = 6M$} & \multicolumn{2}{c||}{Case 1} & \multicolumn{2}{c|}{Case 2} \\
\hline
$K/S_{0}$ & $\sigma_{imp,LVB}$ & $\epsilon_{LVF}$ & $\epsilon_{SLVB}$ & $\epsilon_{SLVF}$ & $\epsilon_{SLVB}$ & $\epsilon_{SLVF}$ \\
\hline
$0.7$ & $14.6017$ & $0.0032$ & $0.0039$ & $0.0042$ & $0.0036$ & $0.0040$ \\
$0.8$ & $11.9199$ & $0.0020$ & $0.0019$ & $0.0021$ & $0.0019$ & $0.0019$ \\
$0.9$ & $10.2664$ & $0.0013$ & $0.0006$ & $0.0013$ & $0.0007$ & $0.0011$ \\
$1$   & $10.8100$ & $0.0011$ & $0.0001$ & $0.0010$ & $0.0006$ & $0.0007$ \\
$1.1$ & $12.6442$ & $0.0011$ & $0.0001$ & $0.0010$ & $0.0013$ & $0.0003$ \\
$1.2$ & $13.9412$ & $0.0012$ & $0.0002$ & $0.0010$ & $0.0003$ & $0.0004$ \\
$1.3$ & $14.8890$ & $0.0012$ & $0.0004$ & $0.0011$ & $0.0008$ & $0.0011$ \\
\hline
\end{tabular}
\end{center}
\begin{center}
\begin{tabular}{|c|c|c||c|c||c|c|}
\hline
\multicolumn{3}{|c||}{$T = 2Y$} & \multicolumn{2}{c||}{Case 3} & \multicolumn{2}{c|}{Case 4} \\
\hline
$K/S_{0}$ & $\sigma_{imp,LVB}$ & $\epsilon_{LVF}$ & $\epsilon_{SLVB}$ & $\epsilon_{SLVF}$ & $\epsilon_{SLVB}$ & $\epsilon_{SLVF}$ \\
\hline
$0.7$ & $10.1742$ & $0.0009$ & $0.0006$ & $0.0010$ & $0.0008$ & $0.0009$ \\
$0.8$ & $ 9.1690$ & $0.0006$ & $0.0003$ & $0.0006$ & $0.0003$ & $0.0004$ \\
$0.9$ & $ 8.9858$ & $0.0004$ & $0.0001$ & $0.0004$ & $0.0002$ & $0.0003$ \\
$1$   & $ 9.6089$ & $0.0004$ & $0.0000$ & $0.0004$ & $0.0001$ & $0.0003$ \\
$1.1$ & $10.6981$ & $0.0004$ & $0.0000$ & $0.0004$ & $0.0000$ & $0.0003$ \\
$1.2$ & $11.6825$ & $0.0004$ & $0.0001$ & $0.0004$ & $0.0000$ & $0.0002$ \\
$1.3$ & $12.4837$ & $0.0004$ & $0.0001$ & $0.0004$ & $0.0000$ & $0.0001$ \\
\hline
\end{tabular}
\end{center}
\caption{Comparison of the approximated implied volatilities $\sigma_{LVB}, \sigma_{LVF}, \sigma_{SLVB}, \sigma_{SLVF}$ for values $m_{1} = 200, m_{2} = 100, \Delta \tau = 1/200, \theta = 1/3, Q =2.$}
\label{Table:ImpliedVolm200N200}
\end{table}
One observes that the size of the $\epsilon_{LVF}$ is similar in Tables~\ref{Table:ImpliedVolm100N200} and~\ref{Table:ImpliedVolm200N200}.
Hence, the experiments indicate that the performance of the calibration procedure is independent of the number of spatial grid points. The difference between the approximations of the fair value by discretizing either \eqref{eq:SemidiscretizationBackwardEquationLV}, \eqref{eq:SemidiscretizationForwardEquationLV}, \eqref{eq:SemidiscretizationBackwardEquation} or \eqref{eq:SemidiscretizationForwardEquation} is always of the size of the temporal discretization error, which is small.

By increasing $m_{1}, m_{2}$ the values $\sigma_{imp,LVB}$ have noticeably changed, which is related to the convergence of $U_{LV,i_{0}}(T)$ to the exact non-discounted fair value $u_{LV}(X_{0},T)$. The differences in implied volatility in Tables~\ref{Table:ImpliedVolm100N200}, \ref{Table:ImpliedVolm100N400}, \ref{Table:ImpliedVolm200N200} reveal that within the LV model and for the current, realistic values of $m_{1}, m_{2}, \tau$ the spatial discretization error is larger than the temporal discretization error (cf.\ Section \ref{TimeDiscretization}). 
Since the calibration procedure from Section \ref{Calibration} matches the fully discrete LV and SLV models up to a difference with the size of the temporal discretization error, one can define an appropriate semidiscretization of the PDE \eqref{eq:BackwardKolmogorovLV}, control the spatial discretization error within the LV model, and then calibrate the SLV model to the LV model such that the fully discrete SLV model matches the market data up to an error which is dominated by the controlled spatial error from semidiscretization within the LV model.

\setcounter{equation}{0}
\section{Conclusion}\label{Conclusion}

In financial practice, SLV models are calibrated to market data for European call and put options by calibrating them to their underlying LV models. Since there is often no closed-form analytical formula available for the fair value of vanilla options under an LV model, the best one can aim for is that the approximations of the fair value given by the two models are identical whenever similar numerical valuation methods are used.
Here, we choose to perform the numerical option valuation by semidiscretizing the respective backward Kolmogorov equations with finite differences.
By making use of an adjoint semidiscretization of the corresponding forward Kolmogorov equations, we derived an expression for the leverage function such that the semidiscretized SLV model is calibrated exactly to the semidiscretized LV model. In order to employ this expression, one has to solve a large non-linear system of ODEs. 
For the actual numerical calibration, temporal discretization of this system by a suitable ADI method is combined with an inner iteration to deal with the non-linearity. 
Our numerical experiments reveal that the fully discrete approximations of the fair value of European call options under the LV and SLV models are always the same up to the size of the temporal discretization error. Since the spatial discretization error is typically much larger than the temporal discretization error, one can control the former one by defining an appropriate semidiscretization of the LV model and then calibrate the fully discrete SLV model de facto exactly to the fully discrete LV model.

\setcounter{equation}{0}
\section*{Acknowledgements}
The work by the first author has been supported financially by a PhD Fellowship of the Research Foundation--Flanders.


\begin{thebibliography}{99}

\bibitem{AP07} \textsc{Andersen, L.~B.~G., Piterbarg, V.~V.} (2007)
Moment explosions in stochastic volatility models.
\textit{Finance Stoch.}, \textbf{11}, 29--50.

\bibitem{A08} \textsc{Andersen, L.~B.~G.} (2008)
Simple and efficient simulation of the Heston stochastic volatility model.
\textit{J. Comp. Finan.}, \textbf{11}, 1--42.

\bibitem{AH11} \textsc{Andreasen, J., Huge, B.} (2011)
Volatility interpolation.
\textit{Risk}, March, 76--79.

\bibitem{C11} \textsc{Clark, I.~J.} (2011)
\textit{Foreign Exchange Option Pricing: A Practitioner's Guide.}
John Wiley \& Sons.

\bibitem{CS88} \textsc{Craig, I.~J.~D., Sneyd, A.~D.} (1988)
An alternating-direction implicit scheme for parabolic equations
with mixed derivatives.
\textit{Comp. Math. Appl.}, \textbf{16}, 341--350.

\bibitem{D94} \textsc{Dupire, B.} (1994)
Pricing with a smile. 
\textit{Risk}, January, 18--20.

\bibitem{ET11} \textsc{Ekstr\"om, E., Tysk, J.} (2011)
Boundary conditions for the single-factor term structure equation.
\textit{Ann. Appl. Prob.}, \textbf{21}, 332--350.

\bibitem{EKO12} \textsc{Engelmann, B., Koster, F., Oeltz, D.} (2012)
Calibration of the Heston stochastic local volatility model: a finite volume scheme. 
\textit{Available at SSRN 1823769}.

\bibitem{G86} \textsc{Gy\"ongy, I.} (1986)
Mimicking the one-dimensional marginal distributions of processes having an Ito differential.
\textit{Probab. Th. Rel. Fields}, \textbf{71}, 501--516.

\bibitem{HH12} \textsc{Haentjens, T., in~'t~Hout, K.~J.} (2012)
Alternating direction implicit finite difference schemes for the Heston--Hull--White partial differential equation.
\textit{J. Comp. Finan.}, \textbf{16}, 83--110.

\bibitem{H09} \textsc{Henry-Labord\`ere, P.} (2009)
Calibration of local stochastic volatility models to market smiles.
\textit{Risk}, September, 112--117.

\bibitem{IHF10} \textsc{in~'t~Hout, K.~J., Foulon, S.} (2010)
ADI finite difference schemes for option pricing in the Heston model
with correlation.
\textit{Int. J. Numer. Anal. Mod.}, \textbf{7}, 303--320.

\bibitem{IHM11} \textsc{in~'t~Hout, K.~J., Mishra, C.} (2011)
Stability of the modified Craig--Sneyd scheme for
two-dimensional convection-diffusion equations with mixed
derivative term.
\textit{Math. Comp. Simul.}, \textbf{81}, 2540--2548.

\bibitem{IHM13} \textsc{in~'t~Hout, K.~J., Mishra, C.} (2013)
Stability of ADI schemes for multidimensional diffusion equations
with mixed derivative terms.
\textit{Appl. Numer. Math.}, \textbf{74}, 83--94.

\bibitem{IHW09} \textsc{in~'t~Hout, K.~J., Welfert, B.~D.} (2009)
Unconditional stability of second-order ADI schemes applied to
multi-dimensional diffusion equations with mixed derivative terms.
\textit{Appl. Numer. Math.}, \textbf{59}, 677--692.

\bibitem{IHW15} \textsc{in~'t~Hout, K.~J., Wyns, M.} (2016)
Convergence of the Modified Craig--Sneyd scheme for two-dimensional convection-diffusion equations with mixed derivative term.
\textit{J. Comp. Appl. Math.}, \textbf{296}, 170--180.

\bibitem{HV03} \textsc{Hundsdorfer, W., Verwer, J.~G.} (2003)
\textit{Numerical Solution of Time-Dependent Advection-Diffusion-Reaction Equations.}
Berlin: Springer.

\bibitem{I15} \textsc{Itkin, A.} (2015)
High-order splitting methods for forward PDEs and PIDEs.
\textit{Int. J. Theor. Appl. Finance}, \textbf{18}, 1550031-1 -- 1550031-24.

\bibitem{L02} \textsc{Lipton, A.} (2002)
The vol smile problem.
\textit{Risk}, February, 61--65.

\bibitem{M14} \textsc{Mishra, C.} (2014)
\textit{Stability of Alternating Direction Implicit Schemes with Application to Financial
Option Pricing Equations}.
PhD thesis, University of Antwerp.

\bibitem{PVF03} \textsc{Pooley, D.~M., Vetzal, K.~R., Forsyth, P.~A.} (2003)
Convergence remedies for non-smooth payoffs in option pricing.
\textit{J. Comp. Finan.}, \textbf{6}, 25--40.

\bibitem{R84} \textsc{Rannacher, R.} (1984)
Finite element solution of diffusion problems with irregular data.
\textit{Numer. Math.}, \textbf{43}, 309--327.

\bibitem{RMQ07} \textsc{Ren, Y., Madan, D., Qian, M.~Q.} (2007)
Calibrating and pricing with embedded local volatility models.
\textit{Risk}, September, 138--143.

\bibitem{R89} \textsc{Risken, H.} (1989)
\textit{The Fokker-Planck Equation: Methods of Solution and Applications,} 2nd edition.
Berlin: Springer.

\bibitem{S87} \textsc{Scott, L.~O.} (1987)
Option pricing when the variance changes randomly: theory, estimation, and an application.
\textit{J. Financ. Quant. Anal.}, \textbf{22}, 419--438.

\bibitem{T11} \textsc{Tachet, R.} (2011)
\textit{Non-parametric model calibration in finance.}
PhD thesis, Ecole Centrale Paris.

\bibitem{TF10} \textsc{Tataru, G., Fisher, T.} (2010)
Stochastic Local Volatility.
\textit{Technical report, Quantitative Development Group, Bloomberg}.

\bibitem{TR00} \textsc{Tavella, D., Randall, C.} (2000)
\textit{Pricing Financial Instruments: The Finite Difference Method.}
New York: John Wiley \& Sons.

\bibitem{VGO14} \textsc{van der Stoep, A.~W., Grzelak, L.~A., Oosterlee, C.~W.} (2014)
The Heston Stochastic-Local Volatility Model: Efficient Monte Carlo Simulation.
\textit{Int. J. Theor. Appl. Finan.}, \textbf{17}, 1450045-1--1450045-30.

\bibitem{W16} \textsc{Wyns, M.} (2016)
Convergence analysis of the Modified Craig--Sneyd scheme for two-dimensional 
convection-diffusion equations with nonsmooth initial data.
To appear in \textit{IMA J. Numer. Anal.}, \url{http://dx.doi.org/10.1093/imanum/drw028}.

\end{thebibliography}
\end{document}